\documentclass[12pt,centertags,oneside]{article}
\usepackage{amsmath,amstext,amsthm,amscd,typearea}
\usepackage{amssymb}
\usepackage{a4wide}
\usepackage[mathscr]{eucal}
\usepackage{mathrsfs}
\usepackage{typearea}
\usepackage{charter}
\usepackage{pdfsync}
\usepackage{url}

\usepackage{xcolor}

\usepackage[a4paper,width=16.2cm,top=3cm,bottom=3cm]{geometry}

\numberwithin{equation}{section}

\newtheorem{theorem}{Theorem}[section]
\newtheorem{definition}[theorem]{Definition}
\newtheorem{proposition}[theorem]{Proposition}
\newtheorem{corollary}[theorem]{Corollary}
\newtheorem{lemma}[theorem]{Lemma}

\newcommand{\supp}{{\rm Supp}}

\newcommand{\vol}{\mathop{\mathrm{vol}}}

\newcommand{\ddc}{dd^c}
\newcommand{\dc}{d^c}

\newcommand{\Reg}{\text{\normalfont Reg}}

\newcommand{\PSH}{{\rm PSH}}

\newcommand{\B}{\mathbb{B}}

\newcommand{\N}{\mathbb{N}}

\newcommand{\R}{\mathbb{R}}

\title{\bf Derivative of volumes of big cohomology classes}
\providecommand{\keywords}[1]{\textbf{\textit{Keywords:}} #1}
\providecommand{\subject}[1]{\textbf{\textit{Mathematics Subject Classification 2010:}} #1}

\author{Duc-Viet Vu}
\newcommand{\Addresses}{{
		\bigskip
		\footnotesize
		\textsc{Duc-Viet Vu, University of Cologne, Division of Mathematics, Department of Mathematics and Computer Science, Weyertal 86-90, 50931, K\"oln.}
		\noindent
		\par\nopagebreak
		\noindent
		\textit{E-mail address}: \texttt{dvu@uni-koeln.de}
}}

\date{\today}
\begin{document}
\maketitle
\begin{abstract} We prove that the partial derivative of the volume function of big classes along any real divisor in a compact K\"ahler manifold is equal to the numerical restricted volume of that class to the divisor. A consequence of our main result is that the divisorial components of the non-K\"ahler locus of a big class lie in fact in the null locus of that class.  
\end{abstract}
\noindent
\keywords {non-pluripolar product}, {restricted volume}, {transcendental Morse inequality}, {big cohomology class}, {plurisubharmonic envelope}.
\\

\noindent
\subject{32U15},  {32Q15}, {14C20}.

%\tableofcontents

%%%%%%%%%%%%%%%%%%%%%%%%%%%%%
%%%%%%%%%%%%%%%%%%%%%%%%%%%%%%%%%%

\section{Introduction}

Let $(X, \omega)$ be a compact K\"ahler manifold of dimension $n$. For every $\alpha \in H^{p,p}(X, \R)$ ($0 \le p \le n$),  we say that $\alpha$ is \emph{pseudoeffective} if $\alpha$ contains a closed positive $(p,p)$-current. A class $\alpha \in H^{1,1}(X, \R)$ is said to be \emph{big} if $\alpha$ contains a K\"ahler current $T$, \emph{i.e.} $T$ is closed positive $(1,1)$-current and $T \ge \delta \omega$ for some constant $\delta>0$. The notion of bigness extends that of an integral class; see \cite{Demailly_analyticmethod,Demailly-Paun}.  %Let $\mathcal{E}, \mathcal{B}, \mathcal{K}$ be the cone of pseudoeffective classes, big classes, or K\"ahler classes in $H^{1,1}(X, \R)$ respectively. The closure $\overline{\mathcal{K}}$ of $\mathcal{K}$ is called the nef cone.
It is also well-known that $\alpha$ is big if and only if the volume $\vol(\alpha)$ is strictly positive. 

Regularity of the volume function on big classes is a topic of interest in complex geometry. It is closely related to other deep properties of big classes or pseudo-effective classes in general; see \cite{Boucksom-Demailly-Paun-Peternell,Boucksom-derivative-volume,Collins-Tosatti-nullloci,ELMNP-amer,ELRNP-survey,WittNystrom-duality} (to cite just a few) for information. The aim of this paper is to relate the  differentiability of volume function with the notion of numerical restricted volume. Let us now recall how to define $\vol(\alpha)$.

%In this paper, by currents we always mean a closed positive current unless otherwise stated.
Let $\alpha \in H^{1,1}(X,\R)$ be pseudoeffective. A closed positive current $T_{\min} \in \alpha$ is said to have \emph{minimal singularities} if for every closed positive current $R \in \alpha$, and for a closed smooth form $\theta \in \alpha$, we can write $R= \ddc u +\theta$, $T_{\min}= \ddc u_{\min}+ \theta$ and there holds $u \le u_{\min}+ C$ on $X$ for some constant $C>0$ ($u,u_{\min}$ are called potentials of $R,T_{\min}$ respectively).  A current with minimal singularities always exists as observed by Demailly: for every closed smooth form $\theta \in \alpha$, let 
$$V_\theta:= \big(\sup \big\{u \in \PSH(X, \theta): u \le 0\big\}\big)^*,$$ 
where $\PSH(X, \theta)$ denotes the set of $\theta$-psh functions on $X$, and $(\cdots)^*$ means the upper semi-continuous regularisation. Thus $V_\theta \in \PSH(X,\theta)$ and $\ddc V_\theta+ \theta$ is a current with minimal singularities. When $\alpha$ is K\"ahler, every current with bounded potentials has minimal singularities. Hence currents with minimal singularities are not unique.

For $(1,1)$-currents $T_1, \ldots, T_m$ on $X$, we denote by $\langle T_1 \wedge \cdots \wedge T_m \rangle$ the non-pluripolar product of $T_1, \ldots, T_m$; see \cite{BEGZ}. By monotonicity of non-pluripolar products (see Theorem \ref{th-main1-monotonicityI} below), the quantity
$$\vol(\alpha):= \int_X \langle T_{\min}^n \rangle$$
 is independent of the choice of a current with minimal singularities $T_{\min} \in \alpha$.   We call $\vol(\alpha)$ \emph{the volume of $\alpha$}. By \cite{Boucksom-volume,Demailly-Ein-Lazarsfeld}, we know that $\alpha$ is big if and only if $\vol(\alpha)>0$, and  if $\alpha= c_1(L)$ the Chern class of a line bundle $L$ on $X$, then 
$$\vol(\alpha)= \limsup_{k \to \infty} \frac{n!}{k^n} \dim H^0(X, L^k)= \lim_{k \to \infty} \frac{n!}{k^n} \dim H^0(X, L^k).$$  
When $\alpha \in H^{1,1}(X, \R)$ is not pseudoeffective, we simply set $\vol(\alpha):= 0$. 
For $1 \le m \le n$ and for every big class $\alpha$, we define 
$\langle \alpha^m \rangle$ to be the class of the current $\langle T_{\min}^m \rangle$ which is again independent of the choice of $T_{\min}$ by the monotonicity of non-pluripolar products.  
In view of the above discussion, it is desirable to extend known properties of the volume function of integral classes to  the cone $\mathcal{E}$ of pseudoeffective classes in $H^{1,1}(X, \R)$. %Our aim is to extend known properties of volume functions for classes in $\mathcal{E} \cap NS_\R(X)$ to $\mathcal{E}$, where  we recall that $NS_\R(X)$ is the Neron-Severi space of $X$ which is the real vector subspace of $H^{1,1}(X,\R)$ generated by Chern classes of line bundles on $X$. 

For every closed positive current $R$, we denote by $\{R\}$ the cohomology class of $R$. If $V$ is an analytic set of pure dimension, we denote by $\{V\}$ the cohomology class of the current of integration along $V$. 
Let $\alpha$ be a big cohomology class. Recall that  the non-K\"ahler locus $E_{nK}(\alpha)$ of $\alpha$ is defined as follows:
$$E_{nK}(\alpha):= \bigcap_{T \in \alpha} E_{+}(T),$$
where the intersection is run over all K\"ahler currents $T \in \alpha$, and $E_{+}(T)$ denotes the set of points at which the Lelong number of $T$ is strictly positive. By \cite{Boucksom_anal-ENS}, we know that $E_{nK}(\alpha)$ is analytic set and there exists a closed positive current $T$ with (almost) analytic singularities in $\alpha$ such that the polar locus of $T$ is equal to $E_{nK}(\alpha)$ (see Definition \ref{def-analysing} below for the notion of almost analytic singularities). 

The non-nef locus $E_{nn}(\alpha)$ is defined to be
$$E_{nn}(\alpha):= \{x \in X: \nu(T_{\min},x)>0\}.$$ 
where $T_{\min}$ is a current with minimal singularities in $\alpha$ and $\nu(T_{\min},x)$ denotes the Lelong number of $T_{\min}$ at $x$. This definition is independent of the choice of $T_{\min}$ (see \cite{Boucksom_anal-ENS}).  Thus $E_{nn}(\alpha)$ is the union of countably many analytic subsets in $X$ by Siu's semi-continuity theorem on Lelong numbers (\cite{Siu}), and $E_{nn}(\alpha) \subset E_{nK}(\alpha)$. Note that
\begin{align}\label{eq-EnnEnk}
E_{nn}(\alpha)= \bigcup_{\epsilon >0} E_{nn}(\alpha+ \epsilon \{\omega\})=  \bigcup_{\epsilon>0} E_{nK}(\alpha+ \epsilon \{\omega\}),
\end{align}
where the unions are taken over every constant $\epsilon>0$; see \cite{Collins-Tosatti-nullloci} or Lemma \ref{le-EnnEnk} below for a proof. If $\alpha$ is the Chern class of some big line bundle $L$, it was known that $E_{nK}(\alpha)= \B_+(L)$ and $E_{nn}(\alpha)= \B_-(L)$; see \cite{Boucksom-these,Tosatti-survey-on-nakamaye} for proofs of these equalities, and  \cite{ELRNP-ann-fourier} for definitions and properties of $\B_+(L)$ and $\B_-(L)$.

Let $V$ be an irreducible analytic subset of dimension $k$ in $X$ and $\alpha$ be a big class. Following \cite{Boucksom-these,Collins-Tosatti-nullloci}, we define \emph{the numerical restricted volume of $\alpha$ to $V$} as follows.  We set 
\begin{align}\label{eq-def-restricted}
\langle \alpha^{k} \rangle|_{X|V}:= \lim_{\epsilon \to 0^+}\int_{\Reg V} \big\langle \big(T_{\min,\alpha+ \epsilon \{\omega\}}|_{\Reg V}\big)^k \big\rangle,
\end{align} 
  where $T_{\min, \alpha+ \epsilon \{\omega\}}$ is a current of minimal singularities in $\alpha+ \epsilon \{\omega\}$, and $\Reg V$ denotes the regular part of $V$. This definition is independent of the choice of currents with minimal singularities (see \cite[Remark 2.2]{Collins-Tosatti-nullloci}). In the formula (\ref{eq-def-restricted}),  we use the following convention: if potentials of $T_{\min,\alpha+ \epsilon\{\omega\}}$ are equal to $-\infty$ on $V$, then we simply define  $T_{\min, \alpha+ \epsilon \{\omega\}}|_{\Reg V}:= 0$. Thus one sees that $$\langle \alpha^k\rangle|_{X|V}=0$$
if $V \subset E_{nn}(\alpha)$ because in this case for every $\epsilon>0$ small enough, there holds $V \subset E_{nn}(\alpha+ \epsilon \{\omega\})$ (by (\ref{eq-EnnEnk})) which is contained in the polar locus of $T_{\min, \alpha+ \epsilon \{\omega\}}$. On the other hand, if $V \not \subset E_{nn}(\alpha)$, then $V \not \subset E_{nK}(\alpha+ \epsilon \{\omega\})$ and  the restrictions of potentials of $T_{\min,\alpha+\epsilon \{\omega\}}$ to $V$ are not identically equal to $-\infty$.

Recall that if $\alpha= c_1(L)$ and $V \not \subset E_{nK}(\alpha)= \B_{+}(L)$, it was proved in \cite{Matsumura-restricted} and \cite{Hisamoto} that $\langle \alpha^{n-1}\rangle|_{X|V}$ is equal to the restricted volume $\vol_{X|V}(L)$ of $L$ to $V$, which is an object of importance in algebraic geometry. If $D= \sum_{j=1}^m \lambda_j D_j$ is a real divisor where $\lambda_j \in \R$ and $D_j$ is an irreducible hypersurface for $1 \le j \le m$, then we put 
$$\langle \alpha^{n-1}\rangle|_{X|D}:= \sum_{j=1}^m \lambda_j \langle \alpha^{n-1} \rangle|_{X|D_j}.$$
Here is our main result.

\begin{theorem}\label{th-main-differentibility} Let $X$ be a compact K\"ahler manifold of dimension $n$.  For every big class $\alpha \in H^{1,1}(X, \R)$ and for every real divisor $D$ and $\gamma:= \{D\}$, there holds
\begin{align}\label{eq-daohamrieng}
\frac{d}{dt}\bigg|_{t=0}\vol(\alpha+ t \gamma)=n \langle \alpha^{n-1}\rangle|_{X|D}.
\end{align}
And if $\supp D \subset E_{nK}(\alpha)$, then
\begin{align}\label{eq-nullloci}
\langle \alpha^{n-1}\rangle|_{X|D}=0.
\end{align}
In other words, every hypersurface in the non-K\"ahler locus of $\alpha$ is contained in the null locus of $\alpha$.
\end{theorem}

%If $D \not \subset E_{nK}(\alpha)$, then $\langle \alpha^{n-1}\rangle|_{X|D}$ is equal to the restricted volume $\vol_{X|D}(\alpha)$ of $\alpha$ to $D$ defined in \cite{WittNystrom-deform}.
We recall that by definition the null locus of $\alpha$ is the union of  irreducible analytic subsets $V$ in $X$ so that 
$$\langle \alpha^{k}\rangle|_{X|V}= 0.$$
One sees from the proof of Theorem \ref{th-main-differentibility} (see (\ref{eq-tinhlainumericalrestricted2}) there) that the function $\langle \alpha^{n-1}\rangle|_{X|D}$ is continuous when $\alpha$ varies on the big cone. We refer to Proposition \ref{pro-themvaoepsilonvanok} at the end of the paper for more information.

Theorem \ref{th-main-differentibility} extends \cite[Theorem C]{WittNystrom-deform} by Witt Nystr\"om to the case where $D$ is not necessarily smooth and more importantly could be contained in the non-K\"ahler locus of $\alpha$. It seems that the smoothness and the fact that the hypersurface is not contained in $E_{nK}(\alpha)$ was used in a crucial way in the proof of \cite[Theorem C]{WittNystrom-deform}.  As commented in \cite{WittNystrom-deform}, one sees that  \cite[Theorem C]{WittNystrom-deform} (hence Theorem \ref{th-main-differentibility}) extends previous results by Boucksom-Favre-Johsson \cite{Boucksom-derivative-volume} and Lazarsfeld-Musta\c{t}\u{a} \cite{Lazarsfeld-Mustata} in the case where $\alpha$ is integral.

Theorem \ref{th-main-differentibility} is a special case of a conjecture due to Boucksom-Demailly-P\u{a}un-Peternell \cite{BDPP} that for any $\gamma \in H^{1,1}(X,\R)$ there holds
\begin{align}\label{eq-BDPP}
\frac{d}{dt}\big|_{t=0}\vol(\alpha+ t \gamma)=n \langle \alpha^{n-1}\rangle \cdot \gamma.
\end{align}

The case where $\alpha, \gamma$ are in the Neron-Severi space $NS_{\R}(X)$ was proved in \cite{Boucksom-derivative-volume} and \cite{BDPP}. If $X$ is projective, then (\ref{eq-BDPP}) was established by Witt Nystr\"om \cite{WittNystrom-duality}. The conjecture (\ref{eq-BDPP}) is equivalent to other conjectural properties of the big cone such as the orthogonal properties, weak transcendental Morse inequality and the duality of the pseudoeffective cone and the movable ones. We refer to the appendix by Boucksom in \cite{WittNystrom-duality} for a proof of these equivalences. There are numerous works around this topic, for example, apart from those cited above, one can consult \cite{Popovici,Popovici2,Tosatti-weakMorse,Tosatti-orthogonality,Xiao-movable-inter,Xiao-weak-morse}.  We notice that in general one has 
$$ \langle \alpha^{n-1}\rangle \cdot \{D\} \ge  \langle \alpha^{n-1}\rangle|_{X|D},$$
see \cite[Section 8]{Collins-Tosatti-nullloci}. %The equality occurs if $D \not \subset E_{nK}(\alpha)$ (CHECK IT AGAIN!).
We hope that  our method can be improved to determine if $\langle \alpha^{n-1}\rangle \cdot \{D\} =  \langle \alpha^{n-1}\rangle|_{X|D}$.

When $\dim X=2$, the non-K\"ahler locus of $\alpha$ only consists of curves (it has no isolated points), hence Theorem \ref{th-main-differentibility} implies a known result that the null locus of a big class in a surface is equal to its non-K\"ahler locus. This fact is also deduced from a recent result by Collins-Tosatti \cite{Collins-Tosatti-nullloci} that for every big class admitting a Zariski decomposition, the null locus is equal to the non-K\"ahler locus (hence in particular it holds for surfaces, see \cite{Zariski} and also \cite{Boucksom_anal-ENS,Matsumura-restricted}).

In general it was conjectured in \cite{Collins-Tosatti-nullloci} that the null locus of a big class $\alpha$ in a compact K\"ahler manifold is equal to the null locus of $\alpha$. If $\alpha$ is a big real divisor and $X$ is projective, this conjecture was proved by Ein-Lazarsfeld-Musta\c{t}\u{a}-Nakamaye-Popa  \cite[Theorem 5.7]{ELMNP-amer}, see also \cite{Boucksom-Cacciola-Lopez}. The case where $\alpha$ is nef and $X$ is a compact K\"ahler manifold was settled by Collins-Tosatti \cite{Collins-Tosatti}; see also  \cite{Nakamaye,Nakamaye-nef} for the integral nef case. We also refer to \cite{Collins-Tosatti-nullloci,Tosatti-survey-on-nakamaye} for a summary and discussions around this question.

% We recover the following continuity property.  

%\begin{corollary} \label{cor-tinhlientuccuavolume}  The continuity of restricted volume functions for hypersurfaces.???\end{corollary}

We have comments on the proof of Theorem \ref{th-main-differentibility}. Our method is quite different from that in \cite{WittNystrom-deform} for which the construction of deformation of normal cones and the partial Legendre transform of psh functions play a key role. For our proof, we express the numerical restricted volume in terms of relative non-pluripolar products of classes by proving that 
$$\langle \alpha^{k} \rangle|_{X|V}= \lim_{\epsilon \to 0^+}\int_X \{\langle (\alpha+ \epsilon \{\omega\})^{k} \dot{\wedge}[V]\rangle \},$$
(Lemma \ref{le-sosanhhaikhainiemnumerestric}); see  Section \ref{sec-productpseu} for the notion in the right-hand side and its properties. This formula is more convenient to work with.

By concavity of volume functions, proving (\ref{eq-daohamrieng}) is equivalent to showing that the left-derivative at $0$ of $f(t):= \vol(\alpha+ t \gamma)$ satisfies:
$$f'(0^+)= n \langle \alpha^{n-1} \rangle|_{X|D}.$$
In order to prove this equality, we first refine a crucial idea of Witt Nystr\"om \cite{WittNystrom-duality} and use monotonicity of non-pluripolar products to obtain 
$$f'(0^+) \le n \langle \alpha^{n-1} \rangle|_{X|D}.$$
The converse inequality is proved using recent notions of products of pseudoeffective classes introduced in \cite{Vu_lelong-bigclass,Viet-generalized-nonpluri}. Some of arguments in this part of the proof are similar to those in \cite{Vu_lelong-bigclass}.

The paper is organized as follows. We recall basic properties of relative non-pluripolar products and the definition of products of pseudoeffective classes from  \cite{Vu_lelong-bigclass} in Section \ref{sec-productpseu}. Other necessary materials about volume functions are presented there as well. In Section \ref{sec-Morse} we establish a version of transcendental Morse inequality refining that in \cite{WittNystrom-duality}. Theorem \ref{th-main-differentibility} is proved in Section \ref{sec-partial-deri}.
\\

%Questions: Higher derivative of volume functions, behaviour near the boundary of the pseudoeffective cone; what would be "volume" of a non-big line bundle or more generally a non-big pseudoeffective classes.

\noindent
\textbf{Acknowledgments.} The author would like to thank Witt Nystr\"om and Valentino Tosatti for enlightening discussions about matters related to this paper. Thanks are also due to George Marinescu for fruitful discussions about the notion of almost analytic singularities. The research of the author is partially funded by the Deutsche Forschungsgemeinschaft (DFG, German Research Foundation)-Projektnummer 500055552 and by the ANR-DFG grant QuaSiDy, grant no ANR-21-CE40-0016.

\section{Products of big cohomology classes}\label{sec-productpseu}

We first recall some basic facts about relative non-pluripolar products. This notion was introduced in \cite{Viet-generalized-nonpluri} as a generalization of the usual non-pluripolar products  given in \cite{BT_fine_87,BEGZ,GZ-weighted}. To simplify the presentation, we only consider the compact setting.

Let $X$ be a compact K\"ahler manifold of dimension $n$.  Let $T_1, \ldots, T_m$ be closed positive $(1,1)$-currents on $X$. Let $T$ be  a closed positive current of bi-degree $(p,p)$ on $X$.  By \cite{Viet-generalized-nonpluri}, we can define \emph{the $T$-relative non-pluripolar product} $\langle \wedge_{j=1}^m T_j \dot{\wedge} T\rangle$  in a way similar to that of  the usual non-pluripolar product.  For readers' convenience, we recall how to do it. 

Write $T_j= \ddc u_j+ \theta_j$, where $\theta_j$ is a smooth form and $u_j$ is a $\theta_j$-psh function. Put 
$$R_k:=\bold{1}_{\cap_{j=1}^m \{u_j >-k\}} \wedge_{j=1}^m (\ddc \max\{u_j,-k\} + \theta_j)\wedge T$$
for $k \in \N$.  By the strong quasi-continuity of bounded psh functions (\cite[Theorems 2.4 and 2.9]{Viet-generalized-nonpluri}), we have 
$$R_k= \bold{1}_{\cap_{j=1}^m \{u_j >-k\}} \wedge_{j=1}^m (\ddc \max\{u_j,-l\} + \theta_j)\wedge T$$
for every $l \ge k \ge 1$. A similar equality also holds if we use local potentials of $T_j$ instead of global ones. We can show that  $R_k$ is positive (see \cite[Lemma 3.2]{Viet-generalized-nonpluri}). 

As in \cite{BEGZ}, since $X$ is K\"ahler, one can check  that $R_k$ is of mass bounded uniformly in $k$ and $(R_k)_k$ admits a limit current which is closed as $k \to \infty$. The last limit is denoted by $\langle \wedge_{j=1}^m T_j \dot{\wedge} T\rangle$.  The  last product is, hence,  a well-defined closed positive current of bi-degree $(m+p,m+p)$; and it is  symmetric with respect to $T_1, \ldots, T_m$ and homogeneous. If $T= [X]$ the current of integration along $X$, then $\langle \wedge_{j=1}^m T_j \dot{\wedge} T\rangle$ is equal to the usual non-pluripolar product $\langle \wedge_{j=1}^m T_j\rangle$ (more generally see Lemma \ref{le-relative-V} below). % we have the following result which follows directly from the definition of relative non-pluripolar products.

For a $(1,1)$-current $P$, recall that \emph{the polar locus} $I_P$ of $P$ is the set of $x\in X$ so that the potentials of $P$ are equal to $-\infty$ at $x$. We note that by definition, (the trace measure of) the current $ \langle  T_1 \wedge \cdots \wedge T_m \dot{\wedge} T\rangle $ has no mass on $\cup_{j=1}^m I_{T_j}$.  Using directly the definition of relative non-pluripolar products, one sees immediately that
\begin{align}\label{eq-giataihypersurface}
\langle [D] \wedge T_2 \wedge \cdots \wedge T_m \dot{\wedge} T\rangle =0,
\end{align}
for every effective divisor $D$, $(1,1)$-currents $T_2,\ldots, T_m$ and every current $T$ on $X$.

\begin{proposition}\label{pro-sublinearnonpluripolar} (\cite[Proposition 3.5]{Viet-generalized-nonpluri}) On  a compact K\"ahler manifold $X$ there hold: 

$(i)$ The product $ \langle  T_1 \wedge \cdots \wedge T_m \dot{\wedge} T\rangle $ is symmetric with respect to $T_1, \ldots, T_m$. 

$(ii)$   Given a positive real number $\lambda$, we have $\langle (\lambda T_1) \wedge T_2 \wedge \cdots \wedge T_m \dot{\wedge} T \rangle = \lambda \langle T_1 \wedge T_2 \wedge \cdots \wedge T_m \dot{\wedge} T\rangle $.

$(iii)$ Given a (locally) complete pluripolar set $A$ such that $T$ has no mass on $A$, then $\langle  T_1 \wedge T_2 \wedge \cdots \wedge T_m \dot{\wedge} T\rangle $ also has no mass on $A$.

$(iv)$  Let $T'_1$ be  a closed positive $(1,1)$-current on $X$ and $T_j,T$ as above. Then 
\begin{align}\label{ine-convexnonpluripoarTjTphayj}
\big\langle  (T_1+T'_1) \wedge \bigwedge_{j=2}^m T_j \dot{\wedge} T\big\rangle \le \langle T_1 \wedge \bigwedge_{j=2}^m T_j \dot{\wedge}T\rangle+ \langle T'_1 \wedge  \bigwedge_{j=2}^m  T_j \dot{\wedge} T\rangle.
\end{align}
The equality occurs if  $T$ has no mass on $I_{T_1} \cup I_{T'_1}$. % or at least one of $T_1, T'_1$ is of locally bounded potentials.  

$(v)$ Let $1 \le l \le m$ be an integer.  Let $T''_j$ be  a closed positive $(1,1)$-current on $X$ and $T_j,T$ as above for $1 \le j \le l$. Assume that $T''_j \ge T_j$ for every $1 \le j \le l$ and $T$ has no mass on $\cup_{j=1}^l I_{T''_j-T_j}$. Then, we have 
$$\langle \bigwedge_{j=1}^l T''_j \wedge \bigwedge_{j=l+1}^m T_j \dot{\wedge} T \rangle \ge \langle \bigwedge_{j=1}^m T_j \dot{\wedge} T \rangle.$$

$(vi)$ Let $A$ be a locally complete pluripolar set. Then we have 
$$\bold{1}_{X \backslash A}\langle  T_1 \wedge T_2 \wedge \cdots \wedge T_m \dot{\wedge} T\rangle= \big\langle  T_1 \wedge T_2 \wedge \cdots \wedge T_m \dot{\wedge}(\bold{1}_{X \backslash A} T)\big\rangle.$$
In particular,  the equality
$$\langle \bigwedge_{j=1}^m T_j \dot{\wedge} T \rangle = \langle \bigwedge_{j=1}^m T_j \dot{\wedge} T' \rangle$$
holds, where $T':= \bold{1}_{X \backslash \cup_{j=1}^m I_{T_j}} T$.
\end{proposition}

The following result give more details on the relation between the relative non-pluripolar product and the non-pluripolar product.

\begin{lemma}  \label{le-relativeandnonrelative} (\cite[Lemma 2.3]{Vu_lelong-bigclass}) Assume that  $T$ is of bi-degree $(1,1)$.  Then we have 
\begin{align}\label{eq-gennonpluriclassi3}
\langle T_1 \wedge \cdots \wedge T_m  \wedge T \rangle=\langle T_1 \wedge \cdots \wedge T_m \dot{\wedge} ( \bold{1}_{X \backslash I_{T}}T) \rangle,
\end{align}
In particular, if $T$ has no mass on $I_T$, then 
$$\langle T_1 \wedge \cdots \wedge T_m  \wedge T \rangle=\langle T_1 \wedge \cdots \wedge T_m \dot{\wedge} T \rangle.$$
\end{lemma}

The next result emphasizes again the role of relative non-pluripolar products. 

\begin{lemma}\label{le-relative-V} Let $V$ be an irreducible analytic set in $X$ and let $[V]$ denote the current of integration along $V$. Let $T_1,\ldots,T_m$ be closed positive $(1,1)$-currents on $X$. Then the following properties hold:

(i) if $V$ is contained in $\cup_{j=1}^m I_{T_j}$, then  $\langle T_1\wedge \cdots \wedge T_m \dot{\wedge} [V] \rangle =0$ and there is $1 \le j_0 \le m$ so that $V \subset I_{T_{j_0}}$,

(ii) if $V$ is not contained in $\cup_{j=1}^m I_{T_j}$ (hence $V \not \subset I_{T_j}$ for every $1\le j \le m$), then 
\begin{align}\label{eq-pushforwardVrelative}
\langle \wedge_{j=1}^m T_j \dot{\wedge} [V]\rangle =i_* \langle T_1|_V \wedge \cdots \wedge T_m|_V\rangle,
\end{align}
where $i: \Reg V \to X$ is the natural inclusion ($\Reg V$ is the regular part of $V$), and $T_j|_V:= \ddc (u_{j}|_{\Reg V})$ if $\ddc u_j = T_j$ (locally).  % We refer to  \cite[Proposition 3.5]{Viet-generalized-nonpluri} for more properties of relative non-pluripolar products. 
\end{lemma}

\proof We prove (i). If there is no $j$ such that $V \subset I_{T_{j}}$, then $I_{T_j} \cap V$ is of zero Lebesgue measure on $\Reg V$ (because the restriction of every local potentials of $T_j$ to $\Reg V$ is again a psh function on $\Reg V$), hence, the union $\cup_{j=1}^m I_{T_j}$ cannot contain $V$. 
The fact that $\langle T_1\wedge \cdots \wedge T_m \dot{\wedge} [V] \rangle =0$ follows directly from Proposition \ref{pro-sublinearnonpluripolar} (vi).

We check Property (ii). Let $u_j$ be a local potential of $T_j$ on an open subset $U$ in $X$. Since $V \not \subset I_{T_j}$, we have $u_j|_V \not \equiv -\infty$. We claim that $u_j|_{\Reg V} \not \equiv -\infty$. Desingularizing $V$ we obtain a surjective holomorphic map $\pi: \tilde{V} \to V$ such that $\tilde{V}$ is a smooth manifold. Hence $u_j \circ \pi$ is psh on $\pi^{-1}(U \cap V)$. Consequently, $u_j \circ \pi$ cannot be identically equal to $-\infty$ on an open set in $\tilde{V}$. Hence $u_j|_{\Reg V} \not \equiv -\infty$ as desired. This is to say that the non-pluripolar product $\langle T_1|_V \wedge \cdots \wedge T_m|V \rangle$ is  a well-defined closed positive current (of finite mass) on $\Reg V$. 
Thus observe now that both sides of (\ref{eq-pushforwardVrelative}) have no mass on the singular locus of $V$ (because of Proposition \ref{pro-sublinearnonpluripolar} (iii) and the fact that $[V]$ has no mass on the singular locus of $V$). Furthermore one sees that both sides of (\ref{eq-pushforwardVrelative}) are equal on the complement of the singular locus of $V$ in $X$ by the definition of relative non-pluripolar products. 
\endproof

For every closed positive current $R$, we denote by $\{R\}$ the cohomology class of $R$. If $V$ is an analytic set of pure dimension, we denote by $\{V\}$ the cohomology class of the current of integration along $V$.
Recall that for closed positive $(1,1)$-currents $P$ and $P'$ on $X$, we say that $P'$ is less singular than $P$ if for every global potential $u$ of $P$ and $u'$ of $P'$, then $u \le u' +O(1)$.  A class $\beta \in H^{p,p}(X, \R)$ is said to be \emph{pseudoeffective} if $\beta$ contains a closed positive $(p,p)$-current. For $\beta_1,\beta_2 \in H^{p,p}(X,\R)$, we write $\beta_1 \le \beta_2$ if $\beta_2- \beta_1$ is pseudoeffective.

\begin{theorem} \label{th-main1-monotonicityI} (\cite[Theorem 1.1]{Viet-generalized-nonpluri}) (Monotonicity I) Let $X$ be  a compact K\"ahler manifold and $T_1,\ldots, T_m,T$ closed positive currents on $X$ such that $T_j$ is of bi-degree $(1,1)$ for $1 \le j \le m$.  Let $T'_j$ be closed positive $(1,1)$-current in the cohomology class of $T_j$ on $X$ such that $T'_j$ is less singular than $T_j$ for $1 \le j \le m$. Then we have 
$$\{\langle T_1 \wedge \cdots \wedge T_m \dot{\wedge} T \rangle \} \le \{\langle T'_1 \wedge \cdots \wedge T'_m \dot{\wedge} T \rangle\}.$$ 
In particular there holds
$$\{\langle T_1 \wedge \cdots \wedge T_m \rangle \} \le \{\langle T'_1 \wedge \cdots \wedge T'_m \rangle\}.$$ 
\end{theorem}

We refer to \cite{BEGZ,Lu-Darvas-DiNezza-mono,WittNystrom-mono} for the case where $m=\dim X$. Let $\alpha$ be a pseudoeffective $(1,1)$-class and let $T_{\min,\alpha}$ be a current with minimal singularities in $\alpha$. We define
$$\vol(\alpha):= \langle \alpha^n \rangle:= \int_X \langle T_{\min,\alpha}^n \rangle$$
which is independent of the choice of $T_{\min,\alpha}$ by Theorem \ref{th-main1-monotonicityI}. For every non-pseudoeffective class $\alpha \in H^{1,1}(X, \R)$, we simply set $\vol(\alpha):=0$.
We collect here  the following well-known properties of $\vol(\alpha)$.

\begin{theorem} \label{the-tinhchatcuavolume}  Let $X$ be a compact K\"ahler manifold. Then the following properties hold: 

(i) (\cite{Boucksom-volume,Demailly-Paun}) A class $\alpha \in H^{1,1}(X,\R)$ contains a K\"ahler current if and only if $\vol(\alpha)>0$. We say that $\alpha$ is big if these equivalent conditions are satisfied.

(ii) (\cite{Boucksom-volume}) The volume function $\vol: \mathcal{E} \to \R$ is continuous and the function $\vol^{1/n}$ is concave on $\mathcal{E}$, where $\mathcal{E}$ is the cone of pseudoeffective classes in $H^{1,1}(X, \R)$. 

(iii) (\cite{Demailly-Ein-Lazarsfeld,Fujita}) If $\alpha$ is the Chern class of a line bundle $L$ on $X$, then 
$$\vol(\alpha)= \lim_{k \to \infty} \frac{n!}{k^n} \dim H^0(X, L^k).$$
\end{theorem}

%We have some comments on Theorem \ref{the-tinhchatcuavolume}.  

Let $\alpha_1,\ldots, \alpha_m$ be big $(1,1)$-classes on $X$. Recall that by using a monotonicity of relative non-pluripolar products (\cite[Theorem 1.1]{Viet-generalized-nonpluri}), we can define the cohomology class $\{\langle \alpha_1 \wedge \ldots \wedge \alpha_m \dot{\wedge} T \rangle \}$ to be the class of the current $\langle \wedge_{j=1}^m T_{j,\min} \dot{\wedge} T \rangle$, where $T_{j,\min}$ is a current with minimal singularities in $\alpha_j$ for $1 \le j \le m$.    When $T$ is the current of integration along $X$, we see that the class $\{\langle \alpha_1\wedge  \cdots \wedge \alpha_m \dot{\wedge} T \rangle\}$ is equal to the product $\langle \alpha_1 \wedge \cdots \wedge \alpha_m \rangle$ defined in  \cite[Definition 1.17]{BEGZ}. For every pseudoeffective $(1,1)$-class $\alpha$, we define $I_\alpha$ to be the polar locus $I_{T_{\min,\alpha}}$, where $T_{\min,\alpha}$ is a current of minimal singularities in $\alpha$. In our applications later, we will use the product $\{\langle \alpha_1 \wedge \ldots \wedge \alpha_m \dot{\wedge} T \rangle \}$ when $T$ is a current of integration along an effective divisor.    We recall the following.

\begin{proposition} \label{pro-productclasss} (\cite[Proposition 4.6]{Viet-generalized-nonpluri}) Let $\alpha_1,\ldots, \alpha_m$ be big $(1,1)$-classes on $X$. Then the following properties hold:

$(i)$ The product  $\{\langle  \bigwedge_{j=1}^m \alpha_j \dot{\wedge} T \rangle\}$ is symmetric and   homogeneous  in $\alpha_1,$ $\ldots, \alpha_m$.  

$(ii)$ Let  $\alpha'_1$ be a big $(1,1)$-class.  Assume that $T$ has no mass on $I_{\alpha_1}\cup I_{\alpha'_1}$.  Then, we have 
$$\{\langle (\alpha_1+ \alpha'_1) \wedge \bigwedge_{j=2}^m \alpha_j \dot{\wedge} T \rangle\} \ge \{\langle  \bigwedge_{j=1}^m \alpha_j \dot{\wedge} T \rangle\}+ \{\langle  \alpha'_1 \wedge \bigwedge_{j=2}^m \alpha_j \dot{\wedge} T \rangle\}.$$

$(iii)$ Let $1 \le l \le m$ be an integer. Let  $\alpha''_1, \ldots, \alpha''_l$ be a pseudoeffective $(1,1)$-class such that $\alpha''_j \ge \alpha_j$ for $1 \le j \le l$.  Assume that $T$ has no mass on $I_{\alpha''_j- \alpha_j}$ for every $1 \le j \le l$. Then, we have 
$$\{\langle \bigwedge_{j=1}^l \alpha''_j \wedge  \bigwedge_{j=l+1}^m \alpha_j \dot{\wedge} T \rangle\} \ge  \{ \langle \bigwedge_{j=1}^m \alpha_j \dot{\wedge} T \rangle\}.$$

$(iv)$ If $T$ has no mass on proper analytic subsets on $X$, then the product   $\{\langle  \bigwedge_{j=1}^m \alpha_j \dot{\wedge} T \rangle\}$ is continuous on the set of $(\alpha_1,\ldots, \alpha_m)$ such that $\alpha_1, \ldots, \alpha_m$ are big.
%$(vi)$ Let $1 \le l \le m$ be an integer. Assume $\gamma:= \{\langle \bigwedge_{j=l+1}^m \alpha_j \wedge T \rangle \}$ is well-defined. Then, we have $\{\langle \bigwedge_{j=1}^m \alpha_j \wedge T \rangle\} =\{ \langle \bigwedge_{j=1}^l \alpha_j \wedge R \rangle\}$ if one of the two sides is well-defined.

$(v)$  If $T$ has no mass on proper analytic subsets on $X$ and  $\alpha_1, \ldots, \alpha_m$ are big and nef, then we have    
$$\{\langle \bigwedge_{j=1}^m \alpha_j \dot{\wedge} T \rangle\}= \bigwedge_{j=1}^m \alpha_j \wedge \{T\}.$$
\end{proposition}

We refer to \cite{Boucksom-Favre-Jonsson,BEGZ} for related statements in the case where $T$ is the current of integration along $X$.  
We now define another notion of products of big classes which will play a key role in our study of volume function. 

\begin{theorem} \label{th-mono-current11} (\cite[Theorem 2.1]{Vu_lelong-bigclass} or \cite[Remark 4.5]{Viet-generalized-nonpluri}) (Monotonicity II) Let $X$ be  a compact K\"ahler manifold and let $T_1,\ldots, T_m,T$ be closed positive $(1,1)$-currents on $X$.  Let $T'_j$ and $T'$ be closed positive $(1,1)$-currents in the cohomology class of $T_j$ and $T$ respectively such that $T'_j$ is less singular than $T_j$ for $1 \le j \le m$ and $T'$ is less singular than $T$. Then we have 
$$\{\langle T_1 \wedge \cdots \wedge T_m \dot{\wedge} T \rangle \} \le \{\langle T'_1 \wedge \cdots \wedge T'_m \dot{\wedge} T' \rangle\}.$$
\end{theorem}

Let $\alpha_l, \ldots,\alpha_m$ be big $(1,1)$-classes in $X$ and $\beta$ be a pseudoeffective $(1,1)$-class in $X$. Let $T_{j, \min}, T_{\min}$ be currents with minimal singularities in  the classes $\alpha_j, \beta$ respectively, where $1 \le j \le m$.  By Theorem \ref{th-mono-current11}, the class 
$$\big\{\langle T_{1,\min} \wedge \cdots \wedge T_{m, \min} \dot{\wedge} T_{\min} \rangle \big \}$$   is  a well-defined pseudoeffective class which is independent of the choice of  $T_{\min}$ and $T_{j, \min}$ for $l \le j \le m$. We denote the last class by 
$$\langle \alpha_1\wedge  \cdots \wedge \alpha_m \dot{\wedge} \beta \rangle.$$
%For simplicity, when $l=1$, we remove the bracket $\{ \quad  \}$ from the last notation. 

%The following result holds for the class $\big \{\langle T_1 \wedge \cdots \wedge T_{l-1} \wedge  \alpha_l\wedge  \cdots \wedge \alpha_m \dot{\wedge} \beta \rangle \big\}$  but to avoid cumbersome notations (while keeping the essence of the statements), we only write it for $l=1$. 
  
\begin{proposition} \label{pro-productclasss2} (\cite[Proposition 2.4]{Vu_lelong-bigclass})
$(i)$ The product  $\langle  \wedge_{j=1}^m \alpha_j \dot{\wedge} \beta \rangle$ is symmetric and   homogeneous  in $\alpha_1,$ $\ldots, \alpha_m$.  

$(ii)$ If $\beta'$ is a pseudo-effective $(1,1)$-class, then 
$$\langle \wedge_{j=1}^m \alpha_j \dot{\wedge} \beta \rangle + \langle \wedge_{j=1}^m \alpha_j \dot{\wedge} \beta' \rangle \le \langle \wedge_{j=1}^m \alpha_j \dot{\wedge} (\beta+\beta') \rangle.$$

$(iii)$ Let $1 \le l \le m$ be an integer. Let  $\alpha''_1, \ldots, \alpha''_l$ be a pseudoeffective $(1,1)$-class such that $\alpha''_j \ge \alpha_j$ for $1 \le j \le l$.  Assume that there is a current with minimal singularities in $\beta$ having no mass on $I_{\alpha''_j- \alpha_j}$ for every $1 \le j \le l$. Then, we have 
$$\langle \wedge_{j=1}^l \alpha''_j \wedge  \wedge_{j=l+1}^m \alpha_j \dot{\wedge} \beta \rangle \ge  \langle \wedge_{j=1}^m \alpha_j \dot{\wedge} \beta \rangle.$$

$(iv)$ If there is a current with minimal singularities in $\beta$ having no mass on proper analytic subsets on $X$, then the product   $\{\langle  \wedge_{j=1}^m \alpha_j \dot{\wedge} \beta \rangle\}$ is continuous on the set of $(\alpha_1,\ldots, \alpha_m)$ such that $\alpha_1, \ldots, \alpha_m$ are big.

$(v)$ We have 
$$\langle  \wedge_{j=1}^m \alpha_j  \wedge \beta \rangle  \le \langle  \wedge_{j=1}^m \alpha_j \dot{\wedge} \beta \rangle$$
and the equality occurs if there is a current with minimal singularities $P$ in $\beta$ such that $P=0$ on $I_P$. 
\end{proposition}

%The following result will be useful later. 

\begin{lemma}\label{le-sosanhpolarlocus} (\cite[Lemma 2.5]{Vu_lelong-bigclass}) Let $\alpha$ be a big class and let $T_{\alpha,\min}$ be a current with minimal singularities in $\alpha$. Then, the current $\bold{1}_{I_{T_{\alpha, \min}}}T_{\alpha,\min}$ is a linear combination of currents of integration along irreducible hypersurfaces of $X$, where for every Borel set $A$, we denote by $\bold{1}_A$ the characteristic function of $A$.
\end{lemma}

Let $\alpha$ be a big cohomology class. Let $T_{\min,\alpha}$ be a current with minimal singularities in $\alpha$. Decompose $T_{\min,\alpha}= T_1+ T_2$, where $T_1$ is a linear combination of currents of integration along hypersurfaces, and $T_2$ has no mass on hypersurfaces. Observe that $T_1$ is independent of the choice of $T_{\min,\alpha}$. We denote $N(\alpha)$ the cohomology class of $T_1$. Such a class was already introduced in \cite{Boucksom_anal-ENS}. By \cite{Boucksom_anal-ENS}, one sees that $N(\alpha- N(\alpha))=0$, and $T_1$ is the only current in $N(\alpha)$ and $\vol(\alpha- N(\alpha))= \vol(\alpha)$. Hence we put $\tilde{N}(\alpha):= T_1$ (the unique closed positive current in $N(\alpha)$).

\begin{corollary}\label{cor-veZalpha} Let $\alpha$ be a big class. Then currents of minimal singularities in $\alpha- N(\alpha)$ has no mass on its polar locus.
\end{corollary}

\proof Let $R$ be a current with minimal singularities in $\alpha-N(\alpha)$. Since $N(\alpha- N(\alpha))=0$, we see that $R$ has no mass on hypersurfaces. Since $\bold{1}_{I_R} R$ is just a linear combination of hyperfaces (Lemma \ref{le-sosanhpolarlocus}), we obtain that $\bold{1}_{I_R} R=0$.
\endproof

\begin{lemma}\label{le-dotwedgevoiNalpah} Let  $\alpha$ be a big class and let $D$ be an effective real divisor whose support is contained in $\supp \tilde{N}(\alpha)$. Then $\{\langle \alpha^{n-1} \dot{\wedge} [D] \rangle\}=0$.
\end{lemma}
\proof
Let $T_{\min}$ be a current with minimal singularities in $\alpha$. Thus the polar locus of $T_{\min}$ contains the support of $\tilde{N}(\alpha)$. Hence $\supp D$ is a subset of the polar locus of $T_{\min}$.  It follows that $\langle T_{\min}^{n-1} \dot{\wedge} [D]\rangle=0$ by Lemma \ref{le-relative-V}. The desired assertion then follows.
\endproof

Let $\varphi$ be a quasi-psh function on $X$ and $\mathcal{I}(\varphi)$ be the multiplier ideal sheaf of $\varphi$ whose stalk at $x \in X$ consists of holomorphic germs $f$ at $x$ such that $|f|^2 e^{-2 \varphi}$ is locally integrable at $x$. It is well-known that $\mathcal{I}(\varphi)$ is coherent (see \cite{Nadel-multiplier} and \cite{Demailly_analyticmethod}). Thus we can associate to $\mathcal{I}(\varphi)$ an analytic subset $V_\varphi$ which is the set where local generators of $\mathcal{I}(\varphi)$ vanish simultaneously.

\begin{definition} \label{def-analysing} Let $\psi$ and $\varphi$  be quasi-psh functions on $X$ and $c>0$ be a constant. We say that $\psi$ has \emph{almost analytic singularities associated to $(c,\varphi)$} if the following two conditions are fulfilled:

(i) If $f_1,\ldots, f_l$ are  generators of $\mathcal{I}(\varphi)$ over an open subset $U \subset X$, then there exists a bounded function $h$ on $U$ such that $h$ is smooth outside $V_\varphi$ and 
$$\psi= c \log \sum_{j=1}^l |f_j|^2 + h,$$

(ii) For every smooth modification $\pi: X' \to X$ such that $\pi^* \mathcal{I}(\varphi)$ is generated by a normal crossing divisor. Then for every $z' \in X'$ there exist local coordinates $(z'_1,\ldots, z'_n)$ and $1 \le l \le n$ and  positive constants $a_j \ge 0$ for $1 \le j \le l$ such that locally 
$$\psi \circ  \pi (z')= \sum_{j=1}^l a_j \log |z_j'| + \text{a smooth function}.$$

%The analytic set $V$, defined locally by $V:= \{f_1= \cdots =f_l=0\}$, is globally well-defined and is called the \emph{singular locus} of $\psi$.
\end{definition}

The condition that $h$ is smooth outside $V_\varphi$ is indeed superfluous because it can be deduced from (ii). 
Let $T$ be a closed positive $(1,1)$-current and $\varphi$ be a quasi-psh function and $c>0$ be a constant. We say that $T$ has \emph{almost analytic singularities associated to $(c,\varphi)$} if potentials of $T$ have almost analytic singularities associated to $(c,\varphi)$. We sometimes simply say $T$ has almost analytic singularities (or even only ``analytic singularities'' for brevity) if we don't want to mention $c,\varphi$ explicitly. We note that the polar locus of a current $T$ with almost analytic sinuglarities is an analytic set.  
A more general version of Definition \ref{def-analysing} was introduced in \cite{Coman-Marinescu-Nguyen2}.

In this paper when we speak of (almost) analytic singularities, we always mean it in the sense of Definition \ref{def-analysing}. Demailly's analytic regularisation of psh functions (\cite{Demailly_analyticmethod} or \cite{Demailly-Paun}) implies that one can approximate a closed positive $(1,1)$-current $T$ by a sequence of closed positive $(1,1)$-current $T_k$ with (almost) analytic singularities associated to $(\frac{1}{2k}, k \varphi)$, where $\varphi$ is a global potential of $T$. 

To end this section, we recall the following known property.

\begin{lemma} \label{le-EnnEnk} Let  $\alpha$ be a big class. Then there holds
\begin{align*}
E_{nn}(\alpha)= \bigcup_{\epsilon >0} E_{nn}(\alpha+ \epsilon \{\omega\})=  \bigcup_{\epsilon>0} E_{nK}(\alpha+ \epsilon \{\omega\}),
\end{align*}
where the unions are taken over every constant $\epsilon>0$.
\end{lemma}

\proof Since $E_{nn}(\alpha+ \epsilon \{\omega\}) \subset E_{nK}(\alpha+ \epsilon \{\omega\})$, we get
$$\bigcup_{\epsilon >0} E_{nn}(\alpha+ \epsilon \{\omega\})\subset \bigcup_{\epsilon>0} E_{nK}(\alpha+ \epsilon \{\omega\}).$$
Let $T_{\min, \alpha+ \epsilon \{\omega\}}$ be a current with minimal singularities in $\alpha+ \epsilon \{\omega\}$. 
Using Demailly's analytic approximation of psh functions (\cite{Demailly_analyticmethod}), we know that for every constant $\epsilon>0$, there exists a K\"ahler current $T_\epsilon \in \alpha+ \epsilon \{\omega\}$ with analytic singularities  so that $\nu(T_\epsilon, \cdot) \le \nu(T_{\min, \alpha}, \cdot)$. It follows that 
$$\bigcup_{\epsilon>0} E_{nK}(\alpha+ \epsilon \{\omega\}) \subset E_{nn}(\alpha).$$  
As to the converse inclusion, we consider a K\"ahler current $R \in \alpha$ with analytic singularities. Let $x_0 \in E_{nn}(\alpha)$. Let $\delta>0$ be a small constant so that $\nu(R, x_0)< \nu(T_{\min, \alpha},x_0)/(2 \delta)$. Since $R$ is K\"ahler, there exists a small enough constant $\epsilon>0$ satisfying that 
$$R':= (1-\delta) T_{\min,\alpha+\epsilon \{\omega\}}+ \delta R- (1- \delta) \epsilon \omega$$
is a closed positive current in $\alpha$. It follows that 
$$\nu(R', x_0) \ge \nu(T_{\min, \alpha}, x_0) \ge 2 \delta \nu(R, x_0).$$
Consequently,  $\nu(T_{\min,\alpha+ \epsilon\{\omega\}}, x_0)>0$. In other words, $x_0 \in E_{nn}(\alpha+ \epsilon\{\omega\})$ for $\epsilon>0$ small enough. This finishes the proof.
\endproof

\section{A version of transcendental Morse inequality} \label{sec-Morse}

%The main result of this section is Proposition \ref{pro-DdivisorMorse2} giving a Morse type inequality. Such an inequality was obtained previously by Witt Nystr\"om 

For every big class $\beta$, the K\"ahler locus  $Amp(\beta)$ of $\beta$ is, by definition, the complement in $X$ of the non-K\"ahler locus $E_{nK}(\beta)$ of $\beta$.  Let $D$ be an effective divisor on $X$. We denote by $[D]$ the current of integration along $D$, and by $\{D\}$ the cohomology class of $[D]$. 

We start with the following variant of a crucial estimate  in \cite{WittNystrom-duality}. We fix a norm $\| \cdot \|$ in $H^{1,1}(X,\R)$.

\begin{lemma}\label{le-DdivisorMorse} Let $X$ be a compact K\"ahler manifold of dimension $n$. Let $\alpha$ be a nef $(1,1)$-class and $D$ an effective real divisor. Then for $t > 0$ small enough we have
\begin{align}\label{ine-infinitesimal}
\vol(\alpha- t \{D\}) - \vol(\alpha) \ge - nt \alpha^{n-1} \cdot \{D\} - M t^2,
\end{align}
where $M>0$ is a uniform constant depending only on $\omega$ and  upper bounds of $\vol(\alpha)$ and  $\vol(D)$.
%$$\vol(\alpha - \{D\}) \ge  \vol(\alpha)- n \alpha^{n-1}  \cdot \{D\}.$$
\end{lemma}

\proof Adding a small K\"ahler form to $\alpha$ and using the continuity of the volume function, we can assume that $\alpha$ is K\"ahler.  
Observe that there exists a constant $M>1$ depending only on the norm of $\{D\}$ (equivalently $\vol(D)$) such that there is a smooth form $\eta \in \{D\}$ satisfying 
$$- M \omega \le \eta \le  M \omega.$$ 
Let $g$ be a negative $\eta$-psh function such that $\ddc g + \eta= [D]$. Note that $g$ is smooth outside the support of $D$. 

Let $r>0$ be a constant. Let $\chi_r$ be convex increasing smooth function regularising the max function $\max\{\cdot, -r\}$ such that $\chi_r(t)=-r $ for $t \le -r -1$, $\chi_r(t)=t$ for $-r+1 \le t \le 0$. We have that  $\chi_r(t)$ decreases to  $t$ as $r \to \infty$ for every $t \in \R_{\le 0}$.  Let  $g_r:= \chi_r(g)$. Observe that 
$$\ddc g_r = \chi''_r(g) d g \wedge \dc g + \chi'_r(g) \ddc g= \chi''_r(g) d g \wedge \dc g + \chi'_r(g) (\ddc g+ \eta)- \chi'_r(g) \eta.$$ 
Thus
\begin{align}\label{ine-grddc}
\ddc g_r+ \chi'_r(g) \eta \ge 0.
\end{align}
We note that $\chi'_r(g)$ increases to $1$ pointwise as $r \to \infty$, and $\chi'_r(g)=1$ on the open set $\{g > -r+1\}$.

Let $0<t_0<1$ be a constant such that $\alpha- t {D}$ is big for every $t \in [0,t_0]$.  Let $t \in [0,t_0]$. Let $\theta \in \alpha$ be a smooth K\"ahler form. Now define 
$$\varphi_r:=P_\theta(t g_r):= \bigg(\sup\{v \in \PSH(X,\theta): v \le t g_r\}\bigg)^*.$$
Observe that since $g_r$ is bounded, $\varphi_r$ is a bounded $\theta$-psh function.
Furthermore the sequence $(\varphi_r)_r$ is decreasing as $r \to \infty$ because $g_r$ is decreasing as $r \to \infty$. Moreover by \cite[Lemma 4.2]{WittNystrom-duality}, the limit $\varphi:= \lim_{r\to \infty} \varphi_r$ is a well-defined $\theta$-psh function and $\varphi -t g$ is a $(\theta-t \eta)$-psh function which is of minimal singularities in $\alpha- t \{D\}$. Since the class of $\theta- t\eta$ is $\alpha- t\{D\}$, it follows that
\begin{align*}
\vol(\alpha- t\{D\}) &= \int_X \langle \big(\ddc (\varphi-tg)+ \theta- t\eta\big)^n \rangle\\
&= \int_{X\backslash D} \langle  \big(\ddc (\varphi-tg)+ \theta- t\eta\big)^n \rangle\\
&= \int_{X \backslash D} \langle  \big(\ddc \varphi+ \theta\big)^n\rangle
\end{align*}
because $\ddc g+ \eta= 0$ outside $D$ (by abuse of notation we also use $D$ to denote the support of the divisor $D$). Let $U$ be a relatively compact open subset of $Amp(\alpha -t_0\{D\}) \backslash D$ (which is a subset of $Amp (\alpha - t \{D\})\backslash D$).  Now since $\varphi_r$ decreases to $\varphi$ which is locally bounded on $Amp(\alpha - t_0\{D\}) \backslash D$, using the continuity of Monge-Amp\`ere operators of bounded potentials, one gets
\begin{align}\label{ine-chanduoivolumeofalphaD}
\vol(\alpha -t\{D\})= \int_{X \backslash D} \langle \big(\ddc \varphi+ \theta\big)^n \rangle  \ge \int_{\overline U}\big(\ddc \varphi+ \theta\big)^n \ge \limsup_{r\to \infty}\int_{\overline U} \big(\ddc \varphi_r+ \theta\big)^n.
\end{align}

By \cite{DiNezzaTrapani-support-envelope,Tosatti-envelop}  (see also \cite{Berman-C11regula,Chu-Zhou-optimal-envelope}), we have that $(\ddc \varphi_r + \theta)^n$ is supported on $E_r:= \{\varphi_r= t g_r\}$ (because $g_r$ is smooth) and 
$$\big(\ddc \varphi_r+ \theta\big)^n= \bold{1}_{E_r} (t\ddc g_r+   \theta)^n.$$
and 
\begin{align}\label{ine-grddctungdiem}
t\ddc g_r + \theta \ge 0
\end{align}
pointwise on $E_r \cap Amp(\alpha - \{D\})$ (see Lemma \ref{le-pointwiseenvelope} below). On the other hand, observe that 
$$t\ddc g_r + \theta \le t\ddc g_r+ \theta + t\chi'_r(g) \eta+ Mt \omega$$
pointwise because 
$$ t\chi'_r(g) \eta+ Mt \omega= t \chi'_r(g)(\eta+ M \omega)+ Mt(1- \chi'_r(g)) \omega \ge  0$$ 
by the choice of $M$. Hence by this and (\ref{ine-grddctungdiem}), we obtain
$$\bold{1}_{\{\varphi_r= t g_r\}}(t\ddc g_r + \theta)^n \le (t\ddc g_r+ \theta + t\chi'_r(g) \eta+ M t \omega)^n.$$  Thus 
\begin{align}\label{ine-uocluongtichphanUngang}
\int_{\overline U} \big(\ddc \varphi_r+ \theta\big)^n &= \int_{X} \big(\ddc \varphi_r+ \theta\big)^n- \int_{X \backslash \overline U} \big(\ddc \varphi_r+ \theta\big)^n\\
\nonumber
&= \langle \alpha^n \rangle - \int_{(X \backslash \overline U)\cap \{\varphi_r= t g_r\}} \big(\ddc \varphi_r+ \theta\big)^n\\
\nonumber
&= \vol(\alpha) - \int_{(X \backslash \overline U)\cap \{\varphi_r= t g_r\}} \big(t\ddc g_r+  \theta\big)^n\\
\nonumber
& \ge \vol(\alpha) - \int_{(X \backslash \overline U)\cap \{\varphi_r= t g_r\}} \big(t\ddc g_r+ t\chi_r'(g)\eta+ \theta + Mt \omega\big)^n.
\end{align}
Let $I_r$ be the second term in the right-hand side of the last inequality. Since $\theta \ge 0$ and $\ddc g_r+ \chi'_r(g) \eta \ge 0$, one gets
\begin{align*}
I_r  &\le \int_{X \backslash \overline U} (\theta+t M \omega)^n+  \sum_{j=1}^{n} \binom{n}{j}\int_{X} (\theta+ M t \omega)^{n-j} \wedge (t\ddc g_r+ t\chi'_r(g) \eta)^j\\
&\le \int_{X \backslash \overline U} (\theta+ M t \omega)^n+ n  \int_{X} (\theta+ M t \omega)^{n-1} \wedge (t\ddc g_r+ t\chi'_r(g) \eta)+\\
&\quad  \sum_{j=2}^{n} \binom{n}{j}\int_{X} (\theta+ M t \omega)^{n-j} \wedge (t\ddc g_r+ t\chi'_r(g) \eta)^j\\
%&\le \int_{X \backslash \overline U} (\theta+ M t \omega)^n+ n  \int_{X} (\theta+ M t \omega)^{n-1} \wedge (t\ddc g_r+ t\chi'_r(g) \eta)+ \\
%&\quad  \sum_{j=2}^{n} \binom{n}{j}\int_{X} (\theta+ M t \omega)^{n-j} \wedge (t\ddc g_r+t M \omega)^j\\
&\le \int_{X \backslash \overline U}  (\theta+t M \omega)^n+ n  \int_{X} (\theta+ M t \omega)^{n-1} \wedge (t\ddc g_r+ t\chi'_r(g) \eta)+\\
&\quad  M^n t^2 \sum_{j=2}^{n} \binom{n}{j}\int_{X} (\theta+ M t \omega)^{n-j} \wedge \omega^j.
\end{align*}
Letting $r \to \infty$ gives
$$\limsup_{r\to\infty} I_r \le  \int_{X \backslash \overline U} (\theta+ M t \omega)^n+ n  \int_{X} (\theta+ M t \omega)^{n-1} \wedge (t\ddc g+ t\eta)+ O(t^2)$$
which is equal to
$$\int_{X \backslash \overline U} (\theta+ M t \omega)^n+ n t  \int_{X} (\theta+ M t \omega)^{n-1} \wedge [D]+ O(t^2).$$
This together with (\ref{ine-uocluongtichphanUngang}) and (\ref{ine-chanduoivolumeofalphaD}) and letting $U$ converges to $Amp(\alpha- \{D\})$ yields
$$\vol(\alpha- t\{D\}) \ge \vol(\alpha)- t \alpha^{n-1} \cdot \{D\} - O(t^2).$$
%This combined with (\ref{ine-chanduoivolumeofalphaD}) gives
%\begin{align}\label{ine-chanduoivolumeofalphaD2}
%\vol(\alpha -t\{D\}) & \ge  \limsup_{r\to \infty}\int_{\overline U} \big(\ddc \varphi_r+ \theta\big)^n.\end{align}
Hence (\ref{ine-infinitesimal}) follows.
\endproof

We recall the following well-known observation (see the proof of \cite[Corollary 1.12]{BoucksomBerman}).

 \begin{lemma} \label{le-pointwiseenvelope} We have $\ddc g_r+ \theta \ge 0$ pointwise  on $E_r \cap Amp(\alpha- \{D\})$.
\end{lemma}

\proof Let $x_0$ be a point in $X$ with $\varphi_r(x_0)= g_r(x_0)$. Note that $\varphi_r$ is locally bounded on $Amp(\alpha- \{D\})$.  

We work locally near $x_0$. Consider a local chart $U$ near $x_0$ so that $\theta= \ddc v$ on $U$. Let $u:= \varphi_r- g_r$. Observe that $u(x_0)=0$ and $u(x) \le 0$ for $x \in U$.  We prove the desired assertion by contradiction.   Suppose that $\ddc g_r(x_0) + \theta$ is not a (semi)positive form. Thus there exists a  complex line passing $Y$ through $x_0$ so that the restriction $\ddc (g_r|_Y)(x_0) + \theta|_Y$ of $\ddc g_r+ \theta$ is strictly negative. Let $u':= u|_{Y}$. By the choice of $x_0$, we see that $u'$ attains its local maximum at $x_0$. On the other hand 
$$\ddc u'(x_0) = (\ddc \varphi_r+ \theta)|_Y - (\ddc g_r + \theta)|_Y > 0.$$
Hence $u'$ is strictly subharmonic on an open neighborhood of $x_0$ in $Y$. It follows that $x_0$ cannot be a local maximum point for $u'$ by the maximum principle, a contradiction. Thus the desired property follows.  
\endproof

%\begin{remark}\label{re-laylopalphatrutD} Let $\alpha$ be a nef in  compact subset $K$ in the big cone. Using Lemmas \ref{le-Radamacher} and \ref{le-DdivisorMorse} we also obtain that 
%$$\vol(\alpha- t \{D\})- \vol(\alpha) \ge - nt \langle (\alpha- t \{D\})^{n-1} \rangle  \cdot \{D\}- M t^2.$$\end{remark}

%Let $\gamma$ be a big class. Let $T$ be a current with minimal singularities in $\gamma$. Decompose $T= [D]+ T'$, where $T'$ has no mass on hypersurfaces. The divisor $D$ is independent of the choice of $T$  and is denoted by $N(\gamma)$ (see \cite{Boucksom_anal-ENS}). 

The following is our main result in this section giving a Morse type inequality. Such an inequality was obtained by Witt Nystr\"om  \cite{WittNystrom-duality} in the case  where $\{D\}$ was assumed to be ample (this assumption was probably crucial for the proof given in  \cite{WittNystrom-duality} or in the first version of this paper in ArXiv).

\begin{proposition}\label{pro-DdivisorMorse2} Let $X$ be a compact K\"ahler manifold. Let $\alpha$ be a big $(1,1)$-class and $D$ a prime divisor, i.e, $D= \lambda V$ for $\lambda \ge 0$ and $V$ an irreducible hypersurface.  Then for $t > 0$ small enough, we have
\begin{align}\label{ine-infinitesimal2}
\vol(\alpha- t \{D\}) - \vol(\alpha)\ge - n t \lim_{\epsilon \to 0^+}\{\langle (\alpha- \epsilon \{\omega\})^{n-1} \dot{\wedge}[D]\rangle\} - M t^2,
\end{align}
where $M$ is a uniform constant depending only on $\omega$ and  upper bounds of $\vol(\alpha)$ and $\vol(D)$. 
%$$\vol(\alpha - \{D\}) \ge  \vol(\alpha)- n \langle \alpha^{n-1}\dot{\wedge} [D]\rangle.$$
\end{proposition}

In the above statement and also below we implicitly identify $H^{n,n}(X, \R)$ with $\R$ through the isomorphism $\beta \longmapsto \int_X \beta$ for $\beta \in H^{n,n}(X, \R)$. 

By Proposition \ref{pro-productclasss} (iii) we note that 
the limit $\lim_{\epsilon \to 0^+}\{\langle (\alpha- \epsilon \{\omega\})^{n-1} \dot{\wedge}[D]\rangle\}$ exists and  \`a priori there holds
$$\lim_{\epsilon \to 0^+}\{\langle (\alpha- \epsilon \{\omega\})^{n-1} \dot{\wedge}[D]\rangle\} \le \{\langle \alpha^{n-1} \dot{\wedge} [D]\rangle\}.$$
Later we see that the equality indeed occurs, see Proposition \ref{pro-themvaoepsilonvanok}.

\proof Let $M$ be the constant in Lemma \ref{le-DdivisorMorse}.  We prove that 
\begin{align}\label{ine-infinitesimal2koepsi}
\vol(\alpha- t \{D\}) - \vol(\alpha)\ge - n t \lim_{\epsilon \to 0^+}\{\langle (\alpha+ \epsilon \{\omega\})^{n-1} \dot{\wedge}[D]\rangle\} - M t^2.
\end{align}
We explain first how to deduce the desired inequality from (\ref{ine-infinitesimal2koepsi}). Indeed applying (\ref{ine-infinitesimal2koepsi}) to $\alpha- \epsilon_0 \{\omega\}$ in place of $\alpha$ for some constant $\epsilon_0>0$, we get 
\begin{align*}
\vol(\alpha- \epsilon_0 \{\omega\}- t \{D\}) - \vol(\alpha- \epsilon_0 \{\omega\})&\ge - n t \lim_{\epsilon \to 0^+}\{\langle (\alpha- \epsilon_0 \{\omega\}+ \epsilon \{\omega\})^{n-1} \dot{\wedge}[D]\rangle\} - M t^2\\
& \ge - n t \{\langle (\alpha- \frac{\epsilon_0}{2} \{\omega\})^{n-1} \dot{\wedge}[D]\rangle\} - M t^2
\end{align*}
by Proposition \ref{pro-productclasss} (iii). Letting $\epsilon_0 \to 0$ and  using the continuity of volume functions on the big cone    gives (\ref{ine-infinitesimal2}).

%If $D$ is contained in $\supp N(\alpha)$ then both sides of (\ref{ine-infinitesimal2koepsi}) are zero for $t>0$ small. Hence ok!
%Assume now $D$ is not contained in $N(\alpha)$. Hence 
%$$\langle (\alpha - N(\alpha))^{n-1} \dot{\wedge} [D] \rangle \le \langle \alpha^{n-1} \dot{\wedge} [D] \rangle.$$
%This combined with the fact that  $\vol(\alpha- N(\alpha)- t\{D\}) \le \vol (\alpha- t\{D\})$ and $\vol(\alpha- N(\alpha))= \vol(\alpha)$ shows that, it suffices to assume that $N(\alpha)=0$ by replacing $\alpha$ by $\alpha- N(\alpha)$. 
We split the proof of (\ref{ine-infinitesimal2koepsi}) into two cases. \\

\noindent
\textbf{Case 1.}  $D$ is not contained in the non-K\"ahler locus $E_{nK}(\alpha)$ of $\alpha$. 

Let $\epsilon >0$ be a small constant. Let $\mu: \widehat X \to X$ be a smooth modification so that 
$$\mu^* \alpha= \widehat \alpha+ \{E\},$$
where $\widehat \alpha$ is K\"ahler, $E$ is an effective divisor on $\widehat X$ so that $\mu(E) \subset E_{nK}(\alpha)$, and 
$$\vol(\widehat \alpha) \ge \vol(\alpha)- \epsilon,$$
and $\mu$ is biholomorphic outside an analytic subset $F$ satisfying $\mu(F)$ is of codimension at least $2$ in $X$ (see \cite{BEGZ}). 
Since $D \not \subset E_{nK}(\alpha)$, we see that 
\begin{align} \label{fact-DEnk}
D \not \subset \mu(E).
\end{align}

Let $\widehat D$  be the strict transform of $D$.  Let $\ddc g+ \eta= [D]$ and $g_r$ be a sequence decreasing to $g$ as in the proof of Lemma \ref{le-DdivisorMorse}. 
Observe that the sequence of currents $\ddc (g_r\circ \mu) + \mu^* \eta$ is of mass uniformly bounded (see, e.g, \cite{Boucksom-these} or \cite{DinhSibony_pullback}). Hence by passing to a subsequence if necessary, we can assume that 
$$T= \lim_{r \to 0} (\ddc g_r\circ \mu + \mu^* \eta)$$
exists. By the construction of $g_r$, we see that $T$ is supported on $\mu^{-1}(\supp D) \cup F$ and $T$ is equal to $[\widehat D]$ outside $F$. It follows that 
$$T= \widehat D+ T',$$
 where $T'$ is a closed positive current, and $\mu(\supp T')$ is contained in an analytic subset of $X$ of codimension at least $2$. Hence we obtain 
\begin{align}\label{eq-tinhcahtDmuphay} 
 \mu_* T' =0, \quad \{\widehat D\} \le \mu^* \{D\}.
 \end{align}

Let $T_{min}$ be a current of minimal singularity in $\alpha$. Observe that $\mu^* T_{min}$ is again a current of minimal singularity in $\mu^* \alpha$.   By abuse of notation we also use $\widehat \alpha$ to denote a smooth K\"ahler form in $\widehat \alpha$. 
Now applying Lemma \ref{le-DdivisorMorse} to $\widehat \alpha$ and $\widehat D$, we get
\begin{align*}
\vol(\widehat \alpha - t\{\widehat D\}) &\ge \vol(\widehat \alpha)- n t\widehat \alpha^{n-1} \cdot \{\widehat D\}- Mt^2\\
&= \vol(\widehat \alpha)- n t \int_{\widehat X} \langle \widehat \alpha^{n-1}\dot{\wedge} [\widehat D] \rangle- Mt^2
\end{align*}

Since $[E]$ is supported on pluripolar set $E$, and $[\widehat D]$ has no mass on $E$ thanks to (\ref{fact-DEnk}), we see that 
\begin{align*}
\int_{\widehat X} \langle \widehat \alpha^{n-1}\dot{\wedge} [\widehat D] \rangle &=  \int_{\widehat X} \langle (\widehat \alpha+ [E])^{n-1}\dot{\wedge} [\widehat D] \rangle\\
& \le  \int_{\widehat X}\langle \mu^* T_{\min}^{n-1}\dot{\wedge} [\widehat D] \rangle\\
& =  \int_{X} \langle T_{\min}^{n-1}\dot{\wedge} [D] \rangle
\end{align*}
by monotonicity of relative non-pluripolar products  (Theorem \ref{th-main1-monotonicityI}). This combined with the fact that $\vol(\widehat \alpha) \ge \vol(\alpha)-\epsilon$ yields
\begin{align*}
\vol(\widehat \alpha - t\{\widehat D\}) &\ge \vol(\widehat \alpha)- n  t \{\langle \alpha^{n-1}\dot{\wedge} [D] \rangle \} -M t^2 \\
&\ge \vol(\alpha)- \epsilon - n t \{\langle \alpha^{n-1}\dot{\wedge} [D] \rangle\} - M t^2.
\end{align*}

We now estimate $\vol(\widehat \alpha - t \{\widehat D\})$. Let $R$ be a current of minimal singularities in $\widehat \alpha - t\{\widehat D\}$. Hence 
$$\vol(\widehat \alpha - t\{\widehat D\}) = \int_{\widehat X} \langle R^n \rangle = \int_{\widehat X} \langle (R+[E])^n \rangle = \int_{X} \langle \big(\mu_*(R+[E])\big)^n \rangle$$
Write $\mu_*(R+E)= \mu_*(R+ E- t T')+ t \mu_* T'$ which is of class 
$$\mu_* (\mu^* \alpha- t\mu^* \{D\})+ t\{\mu_* T'\}= \alpha- t \{D\}$$
 because $\mu_*T'=0$. It follows that 
$$\vol(\widehat \alpha - t\{\widehat D\}) \le \vol(\alpha- t\{D\}).$$
%  Observe that the cohomology class of $\mu_*(R+[E])$ is  in $\mu^* $
Hence 
$$\vol(\alpha- t\{D\}) \ge \vol(\alpha)- \epsilon - n t \{\langle \alpha^{n-1} \dot{\wedge} [D] \rangle\} - M t^2$$
for every $\epsilon>0$. Letting $\epsilon \to 0$ gives the desired inequality. 
\\

\noindent
\textbf{Case 2.} $D$ is contained in $E_{nK}(\alpha)$.

Let $P$ be a K\"ahler current in $\alpha$ so that the polar locus of $P$ is equal to $E_{nK}(\alpha)$ (by \cite{Boucksom_anal-ENS}).  Write $P= \lambda [D]+ P'$, where $\lambda>0$ and $P'$ is a closed positive current whose trace measure  has no mass on $D$. 

If $\supp D \subset \supp \tilde{N}(\alpha)$, then  the left-hand side of (\ref{ine-infinitesimal2koepsi}) is zero for $t>0$ small enough thanks to the property that $\vol(\alpha- N(\alpha))= \vol(\alpha)$. 

We consider now $\supp D \not \subset \supp \tilde{N}(\alpha)$. In other words, the currents of minimal singularities in $\alpha$ have zero Lelong number along $D$. Hence the generic Lelong number of the K\"ahler current $\delta P+ (1- \delta) T_{min}- \lambda \delta [D]$ along $D$ is equal to $0$ for $\delta \in (0,1]$. We infer that the non-K\"ahler locus of the  class $\alpha- \delta \lambda \{D\}$ does not contain $D$ for every $\delta \in (0,1]$.  Applying Case 1 to $\alpha - \delta \lambda \{D\}$ in place of $\alpha$ we obtain
\begin{align*}
\vol(\alpha- \delta\lambda \{D\}- t \{D\}) - \vol(\alpha- \delta \lambda \{D\})&\ge - n t \{\big\langle (\alpha- \delta \lambda \{D\})^{n-1} \dot{\wedge}[D] \big\rangle\} - M t^2\\
& \ge  - n t \{\big\langle (\alpha + \delta A \{\omega\})^{n-1} \dot{\wedge}[D] \big\rangle \} - M t^2,
\end{align*}
by Proposition \ref{pro-productclasss} (iii),  where $A>0$ be a big constant so that $A \{\omega\}+ \lambda \{D\}$ is a K\"ahler class. Letting $\delta  \to 0$ gives
\begin{align*}
\vol(\alpha- t \{D\}) - \vol(\alpha) &\ge - n t \lim_{\epsilon \to 0^+} \{\langle (\alpha+ \epsilon \{\omega\})^{n-1} \dot{\wedge}[D]\rangle\} - M t^2
\end{align*}
as desired. The proof is finished.
\endproof

\section{Partial derivative of the volume} \label{sec-partial-deri}

We begin with the following crucial auxiliary lemma which is a variant of \cite[Lemma 3.1]{BoucksomBermanWitt}. 

\begin{lemma}\label{le-concave} Let $f_k$ be a sequence of concave functions in $\R$ and $g$ be a function on $\R$ such that

(i) $\liminf_{k\to \infty} f_k \ge g$,

(ii) $\lim_{k \to \infty} f_k(0) = g(0)$,

(iii) the left-derivative of $g$ at $0$ exists. 

\noindent
Then $\limsup_{k \to \infty} f'_k(0^-) \le g'(0^-)$.
\end{lemma}

\proof By concavity one has
$$f_k(0)+ f'_k(0^-)t \ge f_k(t)$$
for every $t \le 0$. Hence 
$$\liminf_{k \to \infty} t f'_k(0^-) \ge g(t)- g(0)$$
for $t \le 0$. 
It follows that 
$$\limsup_{k \to \infty} f'_k(0^-) \le \liminf_{t \to 0^-}(g(t)- g(0))/t = g'(0^-).$$ 
\endproof

\begin{lemma}\label{le-derivative-alongintegral}
Let $\alpha$ be a big class and $D$ be a prime divisor. Then we have 
\begin{align}\label{eq-daohamtheoDlimepsi}
\lim_{t \to 0^+}\frac{\vol(\alpha - t\{D\})- \vol(\alpha)}{t}\le  - n  \lim_{\epsilon \to 0^+} \{\langle (\alpha- \epsilon \{\omega\})^{n-1} \dot{\wedge}[D]\rangle\},
\end{align}
\end{lemma}

\proof   Let $\gamma:= \{D\}$.  We will prove a stronger property that
\begin{align}\label{ine-conversealphat}
\lim_{t \to 0^+}\frac{\vol(\alpha- t\gamma) - \vol(\alpha)}{t} \le  - n \lim_{\epsilon \to 0^+}\langle (\alpha- \epsilon \{\omega\})^{n-1} \dot{\wedge}\gamma\rangle.
\end{align}
 If $D \subset \supp \tilde{N}(\alpha)$, then 
$$\vol(\alpha)=\vol(\alpha- N(\alpha)) \le \vol(\alpha- t \gamma) \le \vol(\alpha),$$
for $t>0$ small enough. %, and $D \subset \supp N(\alpha- \epsilon \{\omega\})$ because $N(\alpha) \subset N(\alpha- \epsilon \{\omega\})$.  
Hence 
$$\vol(\alpha- t \gamma)= \vol(\alpha)$$
 for $t>0$ small enough and 
 $$\langle \alpha^{n-1} \dot{\wedge} [D]\rangle=0$$
  (Lemma \ref{le-dotwedgevoiNalpah}). Thus 
$$\lim_{t \to 0^+}\frac{\vol(\alpha - t\gamma)- \vol(\alpha)}{t}= - n \{\langle \alpha^{n-1} \dot{\wedge} [D]\rangle\}.$$
On the other hand since $ \supp D \subset \supp \tilde{N}(\alpha)$, we see that $\delta [D] \le \tilde{N}(\alpha)$ for some small enough constant $\delta>0$. It follows that  $[D]$ is the only closed positive current in $\gamma$. Indeed if $S$ is a current in $\gamma$, then $S_1:= \tilde{N}(\alpha)- \delta [D]+ \delta S$ is a current in $N(\alpha)$, thus, $S_1= \tilde{N}(\alpha)$ which implies that $S=[D]$. Consequently, we infer that $\{\langle \alpha^{n-1} \dot{\wedge} [D]\rangle\}= \langle \alpha^{n-1} \dot{\wedge} \gamma\rangle$. The desired equality (\ref{ine-conversealphat}) follows in this case.

It remains to check the desired inequality (\ref{ine-conversealphat}) in the case where $D \not \subset \supp \tilde{N}(\alpha)$. Let $\alpha':= \alpha- N(\alpha)$. We have $N(\alpha')=0$,  $\vol(\alpha')= \vol(\alpha)$ and 
$$\vol(\alpha' - t \gamma)= \vol\big(\alpha- t \gamma- N(\alpha- t\gamma)+ N(\alpha- t \gamma)-  N(\alpha)\big)$$
which is 
$$\ge \vol\big(\alpha- t \gamma- N(\alpha- t\gamma) \big)= \vol(\alpha- t \gamma)$$
because $N(\alpha- t \gamma) \ge N(\alpha)$ (we used here the fact that $D \not \subset \supp N(\alpha)$ and $\supp D$ is an irreducible analytic set). Hence
\begin{align}\label{ine-sosanhalphaphayalpaha}
\vol(\alpha' - t \gamma) - \vol(\alpha') \ge \vol(\alpha- t \gamma) - \vol(\alpha).
\end{align}
On the other hand let $T$ be a current with minimal singularities in $\alpha- \epsilon \{\omega\}$.  Write
$$T= \tilde{N}(\alpha- \epsilon \omega)+ T'= \tilde{N}(\alpha)+ T''$$
where
$$T'':= T'+ \tilde{N}(\alpha- \epsilon \omega)- \tilde{N}(\alpha)$$
which is a closed positive current in $\alpha' - \epsilon \omega$. Let $R$ be a current with minimal singularities in $\gamma$.  By definition we have 
$$\langle (\alpha- \epsilon \{\omega\})^{n-1} \dot{\wedge}\gamma\rangle = \{\langle T^{n-1} \dot{\wedge} R \rangle \}.$$
By Proposition \ref{pro-sublinearnonpluripolar} (iv), we get
$$\langle T^{n-1} \dot{\wedge} R \rangle \le \langle \tilde{N}(\alpha)^{n-1} \dot{\wedge} R \rangle+ \langle (T'')^{n-1} \dot{\wedge} R \rangle$$
which is equal to 
$$\langle (T'')^{n-1} \dot{\wedge} R \rangle$$
because of (\ref{eq-giataihypersurface}).  We infer
$$\{\langle T^{n-1} \dot{\wedge} R \rangle\} \le \{\langle (T'')^{n-1} \dot{\wedge} R \rangle\} \le \langle (\alpha'- \epsilon \{\omega\})^{n-1} \dot{\wedge} \gamma \rangle$$
by monotonicity (Theorem \ref{th-mono-current11}). In other words
$$\langle (\alpha- \epsilon \{\omega\})^{n-1} \dot{\wedge}\gamma\rangle \le  \langle (\alpha'- \epsilon \{\omega\})^{n-1} \dot{\wedge}\gamma\rangle.$$
In view of this and (\ref{ine-sosanhalphaphayalpaha}), it is enough to prove the desired inequality (\ref{ine-conversealphat}) for $\alpha'$ in place of $\alpha$. In other words we can assume that $N(\alpha)=0$. 
It follows that the currents with minimal singularities in $\alpha$ has no mass on their polar locus (Corollary \ref{cor-veZalpha}). This combined with Proposition \ref{pro-productclasss2} (v) yields the following crucial equality
\begin{align}\label{eq-chuyentukochamsangcham}
\langle \alpha^n\rangle = \langle \alpha^{n-1} \dot{\wedge} \alpha\rangle.
\end{align}

Fix $\epsilon >0$ a constant. Consider $t>0$ so that $\epsilon \{\omega\}- t \gamma$ is a K\"ahler class. Thus by Proposition \ref{pro-productclasss} (iii), there holds
$$\langle \alpha^k \wedge (\alpha- t \gamma)^{n-k} \dot{\wedge} \gamma \rangle \ge \langle (\alpha- \epsilon \{\omega\})^{n-1}\dot{\wedge} \gamma \rangle.$$
Combining this with (\ref{eq-chuyentukochamsangcham}) gives 
\begin{align*}
\langle \alpha^n\rangle &= \langle \alpha^{n-1} \dot{\wedge} \alpha\rangle \ge \langle \alpha^{n-1}\wedge (\alpha- t \gamma)\rangle+ t \langle \alpha^{n-1} \dot{\wedge}  \gamma\rangle \\ \quad &(\text{by  Proposition \ref{pro-productclasss2} (ii)}) \\
& = \langle \alpha^{n-2}\wedge (\alpha- t \gamma) \dot{\wedge} \alpha\rangle+ t \langle \alpha^{n-1} \dot{\wedge}  \gamma\rangle \\
& \ge  \langle \alpha^{n-2}\wedge (\alpha- t \gamma) \dot{\wedge} (\alpha- t \gamma)\rangle+ t \langle \alpha^{n-2}\wedge (\alpha- t \gamma) \dot{\wedge}   \gamma\rangle+ \\
& \quad    t \langle \alpha^{n-1} \dot{\wedge}  \gamma\rangle \\
& \ge \langle \alpha^{n-2}\wedge (\alpha- t \gamma)^2 \rangle+ t \langle \alpha^{n-2}\wedge (\alpha- t \gamma) \dot{\wedge}   \gamma\rangle+ \\
& \quad    t \langle \alpha^{n-1} \dot{\wedge}  \gamma\rangle \\
& \ge \langle \alpha^{n-2}\wedge (\alpha- t \gamma)^2 \rangle+ 2 t \langle(\alpha- \epsilon\{\omega\})^{n-1} \dot{\wedge}   \gamma\rangle .
\end{align*}
Repeating this argument for $\langle \alpha^{n-2}\wedge (\alpha- t \gamma)^2 \rangle$ gives
\begin{align*}
\langle \alpha^n\rangle \ge \vol(\alpha- t \gamma)+ nt \langle(\alpha- \epsilon\{\omega\})^{n-1} \dot{\wedge}   \gamma\rangle.
\end{align*}
We deduce that
$$\limsup_{t \to 0^+} \frac{\vol(\alpha- t \gamma)- \vol(\alpha)}{t} \le -n \langle(\alpha- \epsilon\{\omega\})^{n-1} \dot{\wedge}   \gamma\rangle$$
for every constant $\epsilon >0$. Letting $\epsilon \to 0^+$ gives
$$\limsup_{t \to 0^+} \frac{\vol(\alpha- t \gamma)- \vol(\alpha)}{t} \le -n \lim_{\epsilon \to 0^+} \langle(\alpha- \epsilon\{\omega\})^{n-1} \dot{\wedge}   \gamma\rangle.$$
This finishes the proof. 
\endproof

\begin{proposition}\label{pro-derivative-alongintegral}
Let $\alpha_0$ be a big cohomology class and $D$ be a prime divisor. Let $f(t):= \vol(\alpha+ t \gamma)$, where $\gamma:= \{D\}$.
 Then  we have 
$$f'(0^-)= n \lim_{\epsilon \to 0^+} \langle(\alpha- \epsilon\{\omega\})^{n-1} \dot{\wedge}   \gamma\rangle$$
and 
\begin{align}\label{eq-bangnhauDDnga}
\lim_{\epsilon \to 0^+} \langle(\alpha- \epsilon\{\omega\})^{n-1} \dot{\wedge}   \gamma\rangle= \lim_{\epsilon \to 0^+} \{\langle(\alpha- \epsilon\{\omega\})^{n-1} \dot{\wedge}   [D]\rangle\}.
\end{align}
\end{proposition}

\proof Let $\gamma:= \{D\}$.
 By Proposition \ref{pro-DdivisorMorse2} and Monotonicity II (Theorem \ref{th-mono-current11}), one has
\begin{align*}
\vol(\alpha- t\gamma) - \vol(\alpha)  &\ge -nt  \lim_{\epsilon \to 0^+}\{\langle (\alpha- \epsilon \{\omega\})^{n-1} \dot{\wedge}[D]\rangle\} - M t^2 \\
& \ge  - nt \lim_{\epsilon \to 0^+}\langle (\alpha- \epsilon \{\omega\})^{n-1} \dot{\wedge}\gamma\rangle- Mt^2.
\end{align*}
Hence 
$$f'(0^-) \le n \lim_{\epsilon \to 0^+} \langle(\alpha- \epsilon\{\omega\})^{n-1} \dot{\wedge}   \gamma\rangle.$$
The converse inequality is Lemma \ref{le-derivative-alongintegral}. The equality (\ref{eq-bangnhauDDnga}) follows from the fact that 
\begin{align*}
 - n  \lim_{\epsilon \to 0^+} \{\langle (\alpha- \epsilon \{\omega\})^{n-1} \dot{\wedge}[D]\rangle\}  &\le \lim_{t \to 0^+}\frac{\vol(\alpha - t\gamma)- \vol(\alpha)}{t}\\
 &\le  - n  \lim_{\epsilon \to 0^+}\langle (\alpha- \epsilon \{\omega\})^{n-1} \dot{\wedge}\gamma\rangle
 \end{align*}
which is 
$$\le - n  \lim_{\epsilon \to 0^+}\{\langle (\alpha- \epsilon \{\omega\})^{n-1} \dot{\wedge}[D]\rangle\}$$
by Monotonicity II again (Theorem \ref{th-mono-current11}).  
\endproof

\begin{lemma} \label{le-Dnamtrongphankhongkahler}  For every constant $\epsilon>0$, we have that the non-K\"ahler locus of $\alpha$ is contained in $I_{\alpha- \epsilon \omega}$, and  every irreducible hypersurface $D$ in the non-K\"ahler locus of $\alpha$  satisfies
\begin{align}\label{eq-bangkoDtrongEnK}
\lim_{\epsilon \to 0^+} \{\langle(\alpha- \epsilon\{\omega\})^{n-1} \dot{\wedge}  [D]\rangle\}= 0.
\end{align}
\end{lemma}

\proof
Let $T_\epsilon$ be a current with minimal singularities in $\alpha- \epsilon \omega$. Demailly's analytic approximation gives $T'_\epsilon \in \alpha - \epsilon \omega$ with analytic singularities so that $T'_\epsilon \ge -\epsilon \omega/ 2$  and
$$I_{T'_\epsilon} \subset I_{T_\epsilon}= I_{\alpha- \epsilon \omega}.$$
Observe that $T'_\epsilon + \epsilon \omega \in \alpha$ is a  K\"ahler current with analytic singularities. Thus the non-K\"ahler locus of $\alpha$ is contained in $I_{T'_\epsilon +\epsilon \omega}= I_{T'_\epsilon}$ which is, in turn, a subset in $I_{\alpha- \epsilon \omega}$.

The desired equality (\ref{eq-bangkoDtrongEnK}) follows immediately from the fact that $D$ lies in the polar locus of $\alpha- \epsilon \{\omega\}$ and Proposition \ref{pro-sublinearnonpluripolar} (vi). The proof is finished.
\endproof

\begin{lemma}\label{le-sosanhhaikhainiemnumerestric} Let $V$ be an irreducible analytic subset of dimension $k$. We have
$$\langle \alpha^k\rangle|_{X|V}= \lim_{\epsilon \to 0^+} \{\langle(\alpha+ \epsilon\{\omega\})^{k} \dot{\wedge}  [V]\rangle\}.$$
\end{lemma}

\proof If $V \subset E_{nn}(\alpha)$, then both sides are zero by definition. If $V\not \subset E_{nn}(\alpha)$, then $V \not \subset E_{nK}(\alpha+ \epsilon \{\omega\})$ for every $\epsilon>0$ (Lemma \ref{le-EnnEnk}). We see that the restrictions to $V$ of potentials of the currents with minimal singularities in $\alpha+ \epsilon \{\omega\}$ are not identically equal to $-\infty$ on $V$.  Consequently the desired equality follows from Lemma \ref{le-relative-V}.    
\endproof

\begin{proof}[End of the proof of Theorem \ref{th-main-differentibility}] We consider first the case where  $D$ is a prime divisor (i.e. $D$ is an irreducible hypersurface).  
Let
 $$G(t):=  \lim_{\epsilon \to 0^+} \{\langle (\alpha+ t \gamma - \epsilon \{\omega\})^{n-1} \dot{\wedge}[D] \rangle\}$$ 
 and
 $$G_{\delta}(t):= \lim_{\epsilon \to 0^+} \{\langle (\alpha+ \delta\{\omega\}+ t \gamma - \epsilon \{\omega\})^{n-1} \dot{\wedge}[D] \rangle\}$$
 and
 $$g(t):= \big(\vol(\alpha+ t \gamma)\big)^{1/n}, \quad f_k(t):= \big(\vol(\alpha +\{\omega\}/k+ t \gamma)\big)^{1/n}. $$
By properties of volume functions (Theorem \ref{the-tinhchatcuavolume}), we have
$$f_k \to g$$
pointwise as $k \to \infty$, and $f_k,g$ concave and $f_k \ge g$ for every $k$. Moreover by Proposition \ref{pro-derivative-alongintegral}, we obtain
$$f'_k(0^-)= \frac{G_{1/k}(0)}{n f_k(0)^{1-1/n}}, \quad g'(0^-)= \frac{G(0)}{ n g(0)^{1-1/n}}.$$
Applying Lemma \ref{le-concave} to $f_k, g$ we deduce that
$$\limsup_{k \to \infty} G_{1/k}(0) \le G(0).$$
On the other hand observe that
$$G_{1/k}(0) \ge \{\langle (\alpha+ \{\omega\}/2k)^{n-1} \dot{\wedge}[D] \rangle \} \ge G(0)$$
for every $k$. Hence 
$$\lim_{k \to 0}G_{1/k}(0)=\lim_{k \to \infty} \{\langle(\alpha+ \{\omega\}/2k)^{n-1} \dot{\wedge}[D] \rangle\}= G(0).$$
In other words, by Lemma \ref{le-sosanhhaikhainiemnumerestric}, we have proven that 
\begin{align}\label{eq-tinhlainumericalrestricted}
\langle \alpha^{n-1}\rangle|_{X|D}= \lim_{\epsilon \to 0^+} \{\langle(\alpha- \epsilon\{\omega\})^{n-1} \dot{\wedge}  [D]\rangle\}.
\end{align}
This combined with Lemma \ref{le-Dnamtrongphankhongkahler} gives $ \langle \alpha^{n-1}\rangle|_{X|D}=0$ if $D \subset E_{nK}(\alpha)$. Thus (\ref{eq-nullloci}) follows.  Moreover, one deduces from (\ref{eq-tinhlainumericalrestricted}) that for every sequence of big classes $(\alpha_j)_{j \in \N}$ converging to $\alpha$ there holds
\begin{align}\label{eq-tinhlainumericalrestricted2}
\langle \alpha^{n-1}\rangle|_{X|D}= \lim_{j \to \infty} \{\langle \alpha_j^{n-1} \dot{\wedge}  [D]\rangle\} = \lim_{j \to \infty} \langle \alpha_j^{n-1}\rangle|_{X|D}.
\end{align}
As a consequence, the left-derivative of $g$ is continuous. This coupled with the concavity of $g$ yields that the derivative $g'(t)$ of $g$ is well-defined and is equal to  
$$\frac{\langle (\alpha+t \gamma)^{n-1}\rangle|_{X|D}}{n g(t)^{1- 1/n}} \cdot$$ 
In particular (\ref{eq-daohamrieng}) follows.

We treat now the general case. Write $D= \sum_{j=1}^m \lambda_j D_j$ for $\lambda_j \in \R$ and $D_j$ irreducible hypersurfaces.  Let $\gamma:= \{D\}$ and $\gamma_j:= \{D_j\}$ for $1 \le j \le m$.  Consider the function
$$h(t_1,\ldots, t_m):=\vol(\alpha+ t_1 \lambda_1 \gamma_1+ \cdots +t_m \lambda_m \gamma_m),$$
for $t_j$ with $|t_j|$ small enough,  and $\rho(t):= (t, \ldots, t)$ ($m$ times). Thus
$$h \circ \rho(t)= \vol(\alpha+ t \gamma). $$
By (\ref{eq-tinhlainumericalrestricted2}), the partial derivative $\partial_{t_j} h$ is continuous on an small open neighborhood $U$ of $0 \in \R^m$. Hence $h \in \mathcal{C}^1(U)$. Consequently 
$$(h \circ \rho)'(0)= \sum_{j=1}^m \partial_{t_j}h \circ \rho(0)= n \lambda_j \sum_{j=1}^m \langle \alpha^{n-1}\rangle|_{X |D_j}= n \langle \alpha^{n-1} \rangle|_{X|D}.$$
The proof is finished.
\end{proof}

\begin{proposition} \label{pro-themvaoepsilonvanok} Let $D$ be an effective real divisor. Then we have 
\begin{align}\label{eq-loanxalopcoho}
\lim_{\epsilon \to 0}\{\langle (\alpha+\epsilon \{\omega\})^{n-1} \dot{\wedge}[D]\rangle\} &=\lim_{\epsilon \to 0}\langle (\alpha+\epsilon \{\omega\})^{n-1} \dot{\wedge}\{D\}\rangle\\
\nonumber
&=\langle \alpha^{n-1} \dot{\wedge} \{D\}\rangle=   \{\langle \alpha^{n-1} \dot{\wedge} [D]\rangle\}.
\end{align}
%Hence if $\supp D \subset E_{nK}(\alpha)$, then
%\begin{align}\label{eq-nulllociphay}
%\langle \alpha^{n-1}\rangle'|_{X|D}=0.
%\end{align}
\end{proposition}

\proof %Observe that (\ref{eq-nulllociphay}) is a direct consequence of (\ref{eq-loanxalopcoho}) and Lemma \ref{le-Dnamtrongphankhongkahler}. It remains to check (\ref{eq-loanxalopcoho}).  
We consider, for the moment, the case where $D$ is prime. Recall that by Proposition \ref{pro-derivative-alongintegral} and monotonicity, we already have
\begin{align}\label{eq-lopcuaDngalaonxa0}
\lim_{\epsilon \to 0^+}\langle (\alpha-\epsilon \{\omega\})^{n-1} \dot{\wedge}\{D\}\rangle &=\lim_{\epsilon \to 0^+}\{\langle (\alpha-\epsilon \{\omega\})^{n-1} \dot{\wedge}[D]\rangle \}\\
\nonumber
&\le \{\langle \alpha^{n-1} \dot{\wedge} [D]\rangle \} \le \langle \alpha^{n-1} \dot{\wedge} \{D\}\rangle.
\end{align}
 We now claim the following equality 
\begin{align}\label{eq-lopcuaDbangnhauloanxa}
\lim_{\epsilon \to 0^+}\langle (\alpha+\epsilon \{\omega\})^{n-1} \dot{\wedge}\{D\}\rangle=\lim_{\epsilon \to 0^+}\langle (\alpha-\epsilon \{\omega\})^{n-1} \dot{\wedge}\{D\}\rangle.
\end{align}
We explain how to obtain the desired inequalities from (\ref{eq-lopcuaDbangnhauloanxa}). Note that the left-hand side of (\ref{eq-lopcuaDbangnhauloanxa}) is greater than or equal to 
$$\langle \alpha^{n-1} \dot{\wedge} \{D\}\rangle$$
by monotonicity.  Combining this with  (\ref{eq-lopcuaDngalaonxa0}) and (\ref{eq-lopcuaDbangnhauloanxa}) gives (\ref{eq-loanxalopcoho}). 

It remains to check (\ref{eq-lopcuaDbangnhauloanxa}) whose proof is almost identical to that of Theorem \ref{th-main-differentibility}. One just needs to  literally replace $[D]$ by $\{D\}$. Let
 $$G(t):=  \lim_{\epsilon \to 0^+} \langle (\alpha+ t \gamma - \epsilon \{\omega\})^{n-1} \dot{\wedge}\{D\} \rangle$$
 and
 $$G_{\delta}(t):= \lim_{\epsilon \to 0^+} \langle (\alpha+ \delta\{\omega\}+ t \gamma - \epsilon \{\omega\})^{n-1} \dot{\wedge}\{D\} \rangle$$
 and
 $$g(t):= \big(\vol(\alpha+ t \gamma)\big)^{1/n}, \quad f_k(t):= \big(\vol(\alpha +\{\omega\}/k+ t \gamma)\big)^{1/n}. $$
By properties of volume functions (Theorem \ref{the-tinhchatcuavolume}), we have
$$f_k \to g$$
pointwise as $k \to \infty$, and $f_k,g$ concave and $f_k \ge g$ for every $k$. Moreover by Proposition \ref{pro-derivative-alongintegral}, we obtain
$$f'_k(0^-)= \frac{G_{1/k}(0)}{n f_k(0)^{1-1/n}}, \quad g'(0^-)= \frac{G(0)}{ n g(0)^{1-1/n}}.$$
Applying Lemma \ref{le-concave} to $f_k, g$ we deduce that
$$\limsup_{k \to \infty} G_{1/k}(0) \le G(0).$$
On the other hand observe that
$$G_{1/k}(0) \ge \langle (\alpha+ \{\omega\}/2k)^{n-1} \dot{\wedge}\{D\} \rangle  \ge G(0)$$
for every $k$. Hence 
$$\lim_{k \to 0}G_{1/k}(0)=\lim_{k \to \infty} \langle(\alpha+ \{\omega\}/2k)^{n-1} \dot{\wedge}\{D\} \rangle= G(0).$$
Thus (\ref{eq-lopcuaDbangnhauloanxa}) follows. %, hence so does (\ref{eq-lopcuaDbangnhauloanxa2}).  %By the proof of Theorem \ref{th-main-differentibility} one also has 
%$$\lim_{\epsilon \to 0}\{\langle (\alpha+\epsilon \{\omega\})^{n-1} \dot{\wedge}[D]\rangle \}= \{\langle \alpha^{n-1} \dot{\wedge} [D] \rangle\}.$$
%\begin{align}\label{eq-lopcuaDbangnhauloanxa2}
%\lim_{\epsilon \to 0}\{\langle (\alpha+\epsilon \{\omega\})^{n-1} \dot{\wedge}[D]\rangle \}&=\lim_{\epsilon \to 0}\langle (\alpha+\epsilon \{\omega\})^{n-1} \dot{\wedge}\{D\}\rangle\\
%\nonumber
%&=\langle \alpha^{n-1} \dot{\wedge} \{D\}\rangle.
%\end{align}
%Hence the desired equality (\ref{eq-loanxalopcoho}) follows in the case where $D$ is prime. 

Consider now the general case. Write 
$$D= \sum_{j=1}^m \lambda_j D_j,$$
where $\lambda_j \ge 0$ and $D_j$ is an irreducible hypersurface for $1 \le j \le m$. Let $\gamma:= \{D\}$ and $\gamma_j:= \lambda_j\{D_j\}$. By multilinearity of relative non-pluripolar products and the first part of the proof, one has
\begin{align}\label{eq-thaytuprimedendivisor}
\lim_{\epsilon \to 0}\{\langle (\alpha+\epsilon \{\omega\})^{n-1} \dot{\wedge}[D]\rangle \}= \{\langle \alpha^{n-1} \dot{\wedge} [D] \rangle\}.
\end{align}
Now by Proposition \ref{pro-productclasss2} (ii) and the first part of the proof, we obtain 
\begin{align*}
\lim_{\epsilon \to 0^+} \langle (\alpha+\epsilon \{\omega\})^{n-1} \dot{\wedge}\gamma\rangle  &\le \sum_{j=1}^m  \lim_{\epsilon \to 0^+} \langle (\alpha+\epsilon \{\omega\})^{n-1} \dot{\wedge}\gamma_j\rangle\\
&=\sum_{j=1}^m \lim_{\epsilon \to 0^+} \lambda_j \{\langle (\alpha+\epsilon \{\omega\})^{n-1} \dot{\wedge}[D_j]\rangle \}\\
&= \lim_{\epsilon \to 0^+}\{\langle (\alpha+\epsilon \{\omega\})^{n-1} \dot{\wedge}[D]\rangle \}=\{\langle \alpha^{n-1} \dot{\wedge} [D] \rangle\}.
\end{align*}
On the other hand the monotonicity and (\ref{eq-thaytuprimedendivisor}) give
\begin{align*}
\lim_{\epsilon \to 0^+} \langle (\alpha+\epsilon \{\omega\})^{n-1} \dot{\wedge}\gamma\rangle &\ge \lim_{\epsilon \to 0^+} \langle (\alpha- \epsilon \{\omega\})^{n-1} \dot{\wedge}\gamma\rangle \\
& \ge \lim_{\epsilon \to 0^+} \langle (\alpha- \epsilon \{\omega\})^{n-1} \dot{\wedge}[D]\rangle= \{\langle \alpha^{n-1} \dot{\wedge} [D] \rangle\}. 
\end{align*}
Thus the equality must occur everywhere in the above inequalities. Hence, the desired equality (\ref{eq-loanxalopcoho}) follows.
\endproof

\bibliography{biblio_family_MA,biblio_Viet_papers,bib-kahlerRicci-flow}

\begin{thebibliography}{10}

\bibitem{BT_fine_87}
{\sc E.~Bedford and B.~A. Taylor}, {\em Fine topology, \v{S}ilov boundary, and
  {$(dd^c)^n$}}, J. Funct. Anal., 72 (1987), pp.~225--251.

\bibitem{BoucksomBerman}
{\sc R.~Berman and S.~Boucksom}, {\em Growth of balls of holomorphic sections
  and energy at equilibrium}, Invent. Math., 181 (2010), pp.~337--394.

\bibitem{BoucksomBermanWitt}
{\sc R.~Berman, S.~Boucksom, and D.~Witt~Nystr{\"o}m}, {\em Fekete points and
  convergence towards equilibrium measures on complex manifolds}, Acta Math.,
  207 (2011), pp.~1--27.

\bibitem{Berman-C11regula}
{\sc R.~J. Berman}, {\em Bergman kernels and equilibrium measures for line
  bundles over projective manifolds}, Amer. J. Math., 131 (2009),
  pp.~1485--1524.

\bibitem{Boucksom-these}
{\sc S.~Boucksom}, {\em C\^ones positifs des vari\'et\'es complexes compactes}.
\newblock \url{http://sebastien.boucksom.perso.math.cnrs.fr/publis/these.pdf},
  2002.
\newblock Ph.D. thesis.

\bibitem{Boucksom-volume}
\leavevmode\vrule height 2pt depth -1.6pt width 23pt, {\em On the volume of a
  line bundle}, Internat. J. Math., 13 (2002), pp.~1043--1063.

\bibitem{Boucksom_anal-ENS}
\leavevmode\vrule height 2pt depth -1.6pt width 23pt, {\em Divisorial {Z}ariski
  decompositions on compact complex manifolds}, Ann. Sci. \'{E}cole Norm. Sup.
  (4), 37 (2004), pp.~45--76.

\bibitem{Boucksom-Cacciola-Lopez}
{\sc S.~Boucksom, S.~Cacciola, and A.~F. Lopez}, {\em Augmented base loci and
  restricted volumes on normal varieties}, Math. Z., 278 (2014), pp.~979--985.

\bibitem{Boucksom-Demailly-Paun-Peternell}
{\sc S.~Boucksom, J.-P. Demailly, M.~P\u{a}un, and T.~Peternell}, {\em The
  pseudo-effective cone of a compact {K}\"{a}hler manifold and varieties of
  negative {K}odaira dimension}, J. Algebraic Geom., 22 (2013), pp.~201--248.

\bibitem{BDPP}
\leavevmode\vrule height 2pt depth -1.6pt width 23pt, {\em The pseudo-effective
  cone of a compact {K}\"{a}hler manifold and varieties of negative {K}odaira
  dimension}, J. Algebraic Geom., 22 (2013), pp.~201--248.

\bibitem{BEGZ}
{\sc S.~Boucksom, P.~Eyssidieux, V.~Guedj, and A.~Zeriahi}, {\em
  Monge-{A}mp\`ere equations in big cohomology classes}, Acta Math., 205
  (2010), pp.~199--262.

\bibitem{Boucksom-Favre-Jonsson}
{\sc S.~Boucksom, C.~Favre, and M.~Jonsson}, {\em Valuations and
  plurisubharmonic singularities}, Publ. Res. Inst. Math. Sci., 44 (2008),
  pp.~449--494.

\bibitem{Boucksom-derivative-volume}
\leavevmode\vrule height 2pt depth -1.6pt width 23pt, {\em Differentiability of
  volumes of divisors and a problem of {T}eissier}, J. Algebraic Geom., 18
  (2009), pp.~279--308.

\bibitem{Chu-Zhou-optimal-envelope}
{\sc J.~Chu and B.~Zhou}, {\em Optimal regularity of plurisubharmonic envelopes
  on compact {H}ermitian manifolds}, Sci. China Math., 62 (2019), pp.~371--380.

\bibitem{Collins-Tosatti}
{\sc T.~C. Collins and V.~Tosatti}, {\em K\"{a}hler currents and null loci},
  Invent. Math., 202 (2015), pp.~1167--1198.

\bibitem{Collins-Tosatti-nullloci}
\leavevmode\vrule height 2pt depth -1.6pt width 23pt, {\em Restricted volumes
  on {K}\"{a}hler manifolds}, Ann. Fac. Sci. Toulouse Math. (6), 31 (2022),
  pp.~907--947.

\bibitem{Coman-Marinescu-Nguyen2}
{\sc D.~Coman, G.~Marinescu, and V.-A. Nguy\^en}.
\newblock Work in progress, 2023.

\bibitem{Lu-Darvas-DiNezza-mono}
{\sc T.~Darvas, E.~Di~Nezza, and C.~H. Lu}, {\em Monotonicity of nonpluripolar
  products and complex {M}onge-{A}mp\`ere equations with prescribed
  singularity}, Anal. PDE, 11 (2018), pp.~2049--2087.

\bibitem{Demailly_analyticmethod}
{\sc J.-P. Demailly}, {\em Analytic methods in algebraic geometry}, vol.~1 of
  Surveys of Modern Mathematics, International Press, Somerville, MA; Higher
  Education Press, Beijing, 2012.

\bibitem{Demailly-Ein-Lazarsfeld}
{\sc J.-P. Demailly, L.~Ein, and R.~Lazarsfeld}, {\em A subadditivity property
  of multiplier ideals}, vol.~48, 2000, pp.~137--156.
\newblock Dedicated to William Fulton on the occasion of his 60th birthday.

\bibitem{Demailly-Paun}
{\sc J.-P. Demailly and M.~Paun}, {\em Numerical characterization of the
  {K}\"{a}hler cone of a compact {K}\"{a}hler manifold}, Ann. of Math. (2), 159
  (2004), pp.~1247--1274.

\bibitem{DiNezzaTrapani-support-envelope}
{\sc E.~Di~Nezza and S.~Trapani}, {\em Monge-{A}mp\`{e}re measures on contact
  sets}, Math. Res. Lett., 28 (2021), pp.~1337--1352.

\bibitem{DinhSibony_pullback}
{\sc T.-C. Dinh and N.~Sibony}, {\em Pull-back of currents by holomorphic
  maps}, Manuscripta Math., 123 (2007), pp.~357--371.

\bibitem{ELRNP-ann-fourier}
{\sc L.~Ein, R.~Lazarsfeld, M.~Musta\c{t}\u{a}, M.~Nakamaye, and M.~Popa}, {\em
  Asymptotic invariants of base loci}, Ann. Inst. Fourier (Grenoble), 56
  (2006), pp.~1701--1734.

\bibitem{ELMNP-amer}
\leavevmode\vrule height 2pt depth -1.6pt width 23pt, {\em Restricted volumes
  and base loci of linear series}, Amer. J. Math., 131 (2009), pp.~607--651.

\bibitem{ELRNP-survey}
{\sc L.~Ein, R.~Lazarsfeld, M.~Musta\c{t}\v{a}, M.~Nakamaye, and M.~Popa}, {\em
  Asymptotic invariants of line bundles}, Pure Appl. Math. Q., 1 (2005),
  pp.~379--403.

\bibitem{Fujita}
{\sc T.~Fujita}, {\em Approximating {Z}ariski decomposition of big line
  bundles}, Kodai Math. J., 17 (1994), pp.~1--3.

\bibitem{GZ-weighted}
{\sc V.~Guedj and A.~Zeriahi}, {\em The weighted {M}onge-{A}mp\`ere energy of
  quasiplurisubharmonic functions}, J. Funct. Anal., 250 (2007), pp.~442--482.

\bibitem{Hisamoto}
{\sc T.~Hisamoto}, {\em Restricted {B}ergman kernel asymptotics}, Trans. Amer.
  Math. Soc., 364 (2012), pp.~3585--3607.

\bibitem{Lazarsfeld-Mustata}
{\sc R.~Lazarsfeld and M.~Musta\c{t}\u{a}}, {\em Convex bodies associated to
  linear series}, Ann. Sci. \'{E}c. Norm. Sup\'{e}r. (4), 42 (2009),
  pp.~783--835.

\bibitem{Matsumura-restricted}
{\sc S.-i. Matsumura}, {\em Restricted volumes and divisorial {Z}ariski
  decompositions}, Amer. J. Math., 135 (2013), pp.~637--662.

\bibitem{Nadel-multiplier}
{\sc A.~M. Nadel}, {\em Multiplier ideal sheaves and {K}\"{a}hler-{E}instein
  metrics of positive scalar curvature}, Ann. of Math. (2), 132 (1990),
  pp.~549--596.

\bibitem{Nakamaye-nef}
{\sc M.~Nakamaye}, {\em Stable base loci of linear series}, Math. Ann., 318
  (2000), pp.~837--847.

\bibitem{Nakamaye}
\leavevmode\vrule height 2pt depth -1.6pt width 23pt, {\em Base loci of linear
  series are numerically determined}, Trans. Amer. Math. Soc., 355 (2003),
  pp.~551--566.

\bibitem{Popovici2}
{\sc D.~Popovici}, {\em Sufficient bigness criterion for differences of two nef
  classes}, Math. Ann., 364 (2016), pp.~649--655.

\bibitem{Popovici}
\leavevmode\vrule height 2pt depth -1.6pt width 23pt, {\em Volume and
  self-intersection of differences of two nef classes}, Ann. Sc. Norm. Super.
  Pisa Cl. Sci. (5), 17 (2017), pp.~1255--1299.

\bibitem{Siu}
{\sc Y.~T. Siu}, {\em Analyticity of sets associated to {L}elong numbers and
  the extension of closed positive currents}, Invent. Math., 27 (1974),
  pp.~53--156.

\bibitem{Tosatti-weakMorse}
{\sc V.~Tosatti}, {\em The {C}alabi-{Y}au theorem and {K}\"{a}hler currents},
  Adv. Theor. Math. Phys., 20 (2016), pp.~381--404.

\bibitem{Tosatti-survey-on-nakamaye}
\leavevmode\vrule height 2pt depth -1.6pt width 23pt, {\em Nakamaye's theorem
  on complex manifolds}, in Algebraic geometry: {S}alt {L}ake {C}ity 2015,
  vol.~97.1 of Proc. Sympos. Pure Math., Amer. Math. Soc., Providence, RI,
  2018, pp.~633--655.

\bibitem{Tosatti-envelop}
\leavevmode\vrule height 2pt depth -1.6pt width 23pt, {\em Regularity of
  envelopes in {K}\"{a}hler classes}, Math. Res. Lett., 25 (2018),
  pp.~281--289.

\bibitem{Tosatti-orthogonality}
\leavevmode\vrule height 2pt depth -1.6pt width 23pt, {\em Orthogonality of
  divisorial {Z}ariski decompositions for classes with volume zero}, Tohoku
  Math. J. (2), 71 (2019), pp.~1--8.

\bibitem{Viet-generalized-nonpluri}
{\sc D.-V. Vu}, {\em Relative non-pluripolar product of currents}, Ann. Global
  Anal. Geom., 60 (2021), pp.~269--311.

\bibitem{Vu_lelong-bigclass}
\leavevmode\vrule height 2pt depth -1.6pt width 23pt, {\em Lelong numbers of
  currents of full mass intersection}, Amer. J. Math., 145 (2023),
  pp.~647--665.

\bibitem{WittNystrom-duality}
{\sc D.~Witt~Nystr\"{o}m}, {\em Duality between the pseudoeffective and the
  movable cone on a projective manifold}, J. Amer. Math. Soc., 32 (2019),
  pp.~675--689.
\newblock With an appendix by S\'{e}bastien Boucksom.

\bibitem{WittNystrom-mono}
\leavevmode\vrule height 2pt depth -1.6pt width 23pt, {\em Monotonicity of
  non-pluripolar {M}onge-{A}mp\`ere masses}, Indiana Univ. Math. J., 68 (2019),
  pp.~579--591.

\bibitem{WittNystrom-deform}
{\sc D.~Witt~Nystr\"om}, {\em Deformations of {K}\"ahler manifolds to normal
  bundles and restricted volumes of big classes}.
\newblock \url{arXiv:2103.03660}, 2021.

\bibitem{Xiao-weak-morse}
{\sc J.~Xiao}, {\em Weak transcendental holomorphic {M}orse inequalities on
  compact {K}\"{a}hler manifolds}, Ann. Inst. Fourier (Grenoble), 65 (2015),
  pp.~1367--1379.

\bibitem{Xiao-movable-inter}
\leavevmode\vrule height 2pt depth -1.6pt width 23pt, {\em Movable intersection
  and bigness criterion}, Univ. Iagel. Acta Math.,  (2018), pp.~53--64.

\bibitem{Zariski}
{\sc O.~Zariski}, {\em The theorem of {R}iemann-{R}och for high multiples of an
  effective divisor on an algebraic surface}, Ann. of Math. (2), 76 (1962),
  pp.~560--615.

\end{thebibliography}
\bibliographystyle{siam}

\bigskip

\noindent
\Addresses
\end{document}